\documentclass[final]{siamltex}

\usepackage{lineno}
\usepackage{algorithmic} 
\usepackage{algorithm}
\usepackage{float}

\usepackage{lipsum}

\usepackage{subfigure}

\usepackage{amssymb}
\usepackage{wrapfig}
\usepackage{graphicx}

\usepackage{amsmath}
\usepackage{lineno}

\renewcommand{\vec}[1]{\ensuremath\boldsymbol{#1}}
\title{The Augmented Fast Marching Method for Level Set Reinitialization\thanks{This 
        work was supported by the University at Buffalo, State University of New York.}}

\author{David Salac\thanks{Mechanical and Aerospace Engineering, 318 Jarvis Hall, University at Buffalo SUNY, Buffalo, NY, 14260-4400 ({\tt davidsal@buffalo.edu}).}}

\begin{document}

\maketitle

\begin{abstract}
	Including derivative information in the modelling of moving interfaces has been proposed as one method to increase the accuracy of numerical schemes with
	minimal additional cost. Here a new level set reinitialization technique using the fast marching method is presented. This augmented fast marching method
	will calculate the signed distance function and up to the second-order derivatives of the signed distance function for arbitrary interfaces. In addition
	to enforcing the condition $\|\nabla\phi\|^2=1$, where $\phi$ is the level set function, the method ensures that $\nabla\left(\|\nabla\phi\|\right)^2=0$
	and $\nabla\nabla\left(\|\nabla\phi\|\right)^2=0$ are also satisfied. 
	Results indicate that for both two- and three-dimensional interfaces the resulting level set and curvature field are smooth even for 
	coarse grids. Convergence results show that using first-order upwind derivatives and the augmented fast marching method result in a second-order
	accurate level set and gradient field and a first-order accurate curvature field.
\end{abstract}

\begin{keywords} 
Gradient Augmented Level Set, Fast marching method, Reinitialization, Level set, Numerical method
\end{keywords}

\begin{AMS}

\end{AMS}

\pagestyle{myheadings}
\thispagestyle{plain}

\section{Introduction}
\label{sec:1.0}

	The fast marching method (FMM) was introduced by Sethian \cite{Sethian1996} as an efficient method to solve general front propagation problems
	where the propagation speed is monotonic. Since its introduction the fast marching method has been successfully utilized in
	seismology \cite{popovici1998three}, photolithography \cite{sethian1996fast}, medical imaging \cite{malladi1998real,malladi1996n,popovici20023}, 
	and as a component in other numerical schemes \cite{salac2011level,adalsteinsson1999fast}. The fast marching method is also
	an extremely efficient way to compute the distance to an interface \cite{Chopp2001}. It is this last application which is the focus of this work.

	Recent attention has been focused on increasing the accuracy of the level set method by including level set gradient information \cite{nave2010gradient}.
	Results for situations where the velocity field does not depend on the current interface 
	show that the accuracy can be increased with minimal additional computational effort. Issues do arise
	when the velocity field depends on the current interface description. Take for example the modelling of vesicles in external fluid flows.
	The vesicle membrane exerts bending forces on the surrounding fluid \cite{salac2011level,helfrich1973,schwalbe2010}. 
	These forces are related to the curvature and the variation of the curvature
	along the vesicle membrane. Mathematically this means that the bending forces are fourth order derivatives of the level set function. The level set function
	should be as smooth as possible for all time and still be able to describe small scale features of the vesicle. This can be accomplished by 
	periodically reinitializing the level set. As of yet a numerical method to accurately reinitialize the level set function and its gradients has not yet
	been developed.
	
	In this work the reinitialization of a arbitrary level set function is considered. 
	In addition to obtaining the signed distance function of an interface the method will allow for the
	accurate calculation of up to second order derivatives of the signed distance function. This will result in the calculation of smooth curvature fields, which 
	aids in the stability of numerical methods depending on this quantity and its derivatives. In Sec. \ref{sec:2.0} the original fast marching method for 
	reinitialization is briefly presented. The augmented fast marching method is shown in Sec. \ref{sec:3.0}. Here both the two- and three-dimensional 
	systems are considered. Two-dimensional convergence studies and results are presented in Sec. \ref{sec:4.0} while the three-dimensional results are shown 
	in Sec. \ref{sec:5.0}.

\section{The Fast Marching Method for Reinitialization}
\label{sec:2.0}

	Consider an interface implicitly defined as the zero of a function, $\Gamma(t)=\left\{ \vec{x}:\phi\left(\vec{x},t\right)=0 \right\}$, moving 
	with a monotonic speed of $F$. The time, $\tau$, at which the interface crosses a point $\vec{x}$ is the solution to the 
	Eikonal equation, $F\|\nabla\tau(\vec{x})\|=1$. If the speed of the front is one (\textit{i.e.} $F=1$) then the time at which the interface will cross a 
	point $\vec{x}$ is the distance from the point to the interface. By setting $\phi(\vec{x})=\tau(\vec{x})$ and solving $\|\nabla\phi(\vec{x})\|=1$ using the fast marching 
	method it is possible to obtain the distance function of the interface. Denoting the region enclosed by the interface as the negative of the distance 
	function results in a signed distance function.	
	In this section the application of the original fast marching method to level set reinitialization is presented. More information can be found in references 
	\cite{Sethian1996,Chopp2001}.

\subsection{Basics of the FMM}
\label{sec:2.1}

	A two-dimensional computational domain with a uniform grid spacing of $h$ has an embedded interface, $\Gamma$. The goal is to calculate the signed distance
	function to the interface without any spurious motion of $\Gamma$. To accomplish this
	the Eikonal equation with $F=1$,
	\begin{equation}
		 \|\nabla \phi\|=1,
		 \label{eq:Eikonal}
	\end{equation}		 
	is solved in the entire computational domain. 	
	In the fast marching method upwind derivatives are used to approximate $\nabla \phi$. This enforces 
	a causality on the propagation of information. Consider two points, $\vec{x}$ and $\vec{y}$, where the point $\vec{x}$ is closer to the 
	interface than point $\vec{y}$. Due to the fact that $F=1$ and the use of upwind derivatives the value at any given point only
	depends on those points closer to the interface. Thus the value $\phi(\vec{y})$ may depend on the value of $\phi(\vec{x})$, but the value of 
	$\phi(\vec{x})$ will \textit{never} depend on the value of $\phi(\vec{y})$. This leads to an ordering of the nodes which needs to be maintained 
	throughout the reinitialization procedure.
	
\begin{figure}[h!]
		\centering	
		\includegraphics[width=0.4\textwidth]{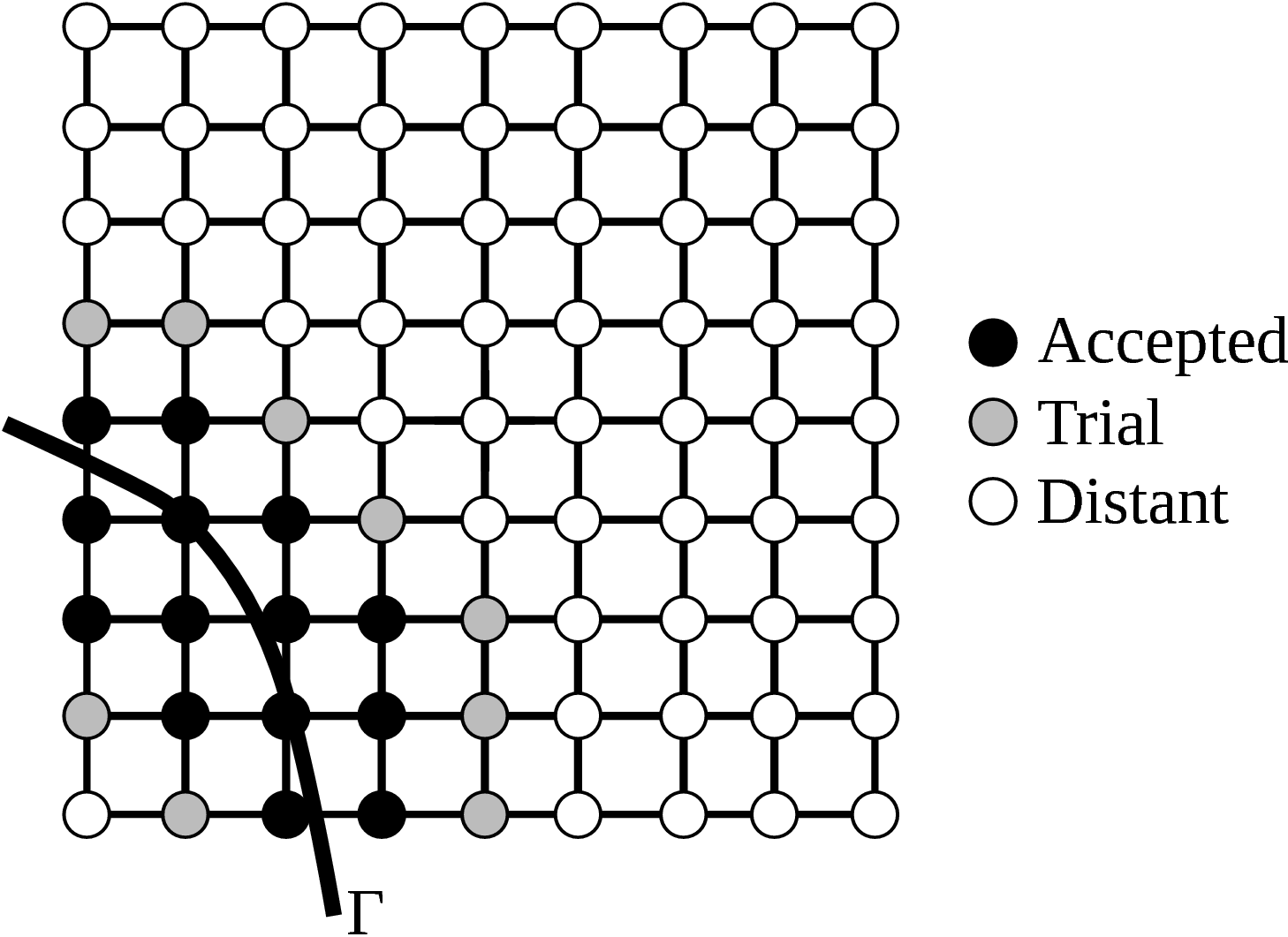}
		\label{fig:2.1.1}
		\caption{The Accepted, Trial, and Distant set of grid nodes.}
	\end{figure}

	When using the fast marching method three sets of nodes are maintained. The first are accepted nodes, $A$, which are those nodes where
	a value of $\phi$ has already been calculated and accepted. The second set are trial nodes, $T$. The trial set contains those nodes which 
	might next be added to the accepted set. Finally, the distant set, $D$, are those nodes which are too far from the interface to be added to the 
	accepted set. Due to the causality of the FMM we are guaranteed that if $\vec{x}\in A$, $\vec{y}\in T$, and $\vec{z}\in D$ 
	then $|\phi(\vec{x})|<|\phi(\vec{y})|<|\phi(\vec{z})|$. See Fig. \ref{fig:2.1.1} for the relationship between the three set of nodes.

\subsection{Initializing the FMM}
\label{sec:2.2}

	To initialize the FMM mark all nodes in the domain as in the distant set. Define a grid cell $\Omega_{i,j}$ as 
	the region enclosed by the four points $\vec{x}_{i,j}$, $\vec{x}_{i+1,j}$, 
	$\vec{x}_{i,j+1}$, and $\vec{x}_{i+1,j+1}$.
	Identify all grid cells which contain the interface and mark the nodes associated with these cells as in the 
	accepted list. 

	The standard fast marching method can be initialized by explicitly calculating the level set function for all nodes associated with cells containing the interface.
	Define a bicubic interpolation function, $P(\vec{x})$, approximating the level set function in a given cell. In the absence of gradient information the bicubic function can be obtained by ensuring 
	that $P(\vec{x}_{m,n})=\phi_{m,n}$, $\partial_x P(\vec{x}_{m,n})=\partial_x \phi_{m,n}$, $\partial_y P(\vec{x}_{m,n})=\partial_y \phi_{m,n}$, and 
	$\partial_{xy} P(\vec{x}_{m,n})=\partial_{xy} \phi_{m,n}$ for $m={i,i+1}$ and $n={j,j+1}$ at grid cell $\Omega_{i,j}$. 
	All derivatives of the level set can be obtained by using standard finite difference approximations. 
	At every grid point $\vec{x}_{i,j}$ in the initially accepted the list the point $\vec{y}$ is calculated such that $p(\vec{y})=0$ and 
	$\nabla p(\vec{y})\times (\vec{x}_{i,j}-\vec{y})=0$.
	The distance to the interface is then $\|\vec{x}-\vec{y}\|_2$. This results in a second-order approximation to the true distance function \cite{Chopp2001}.

\subsection{Updating Nodes in the FMM}
\label{sec:2.3}

	After the initial accepted list is determined the remaining grid nodes are updated in an ordered manner. 
	All grid nodes adjacent to the initial accepted list are given estimates of the distance function by solving an upwind discretization of Eq. (\ref{eq:Eikonal}),
	\begin{equation}	
		\left(D_x^{\pm}\phi\right)^2+\left(D_y^{\pm}\phi\right)^2=1,
		\label{eq:StdFMM}		
	\end{equation}
	where $D_x^{\pm}$ and $D_y^{\pm}$ are one-sided derivatives.
	The appropriate derivative is chosen based on the direction of neighboring accepted nodes. 
	Let the node being updated be $\vec{x}_{i,j}$ with nodes $\vec{x}_{i-1,j}$ and 
	$\vec{x}_{i,j+1}$ in the accepted list. To first order the derivatives in this case would be $D_x^{-}\phi=(\phi_{i,j}-\phi_{i-1,j})/h$ 
	and $D_y^{+}\phi=(\phi_{i,j+1}-\phi_{i,j})/h$, where $h$ is the grid spacing. In general solving Eq. (\ref{eq:StdFMM}) will result in two real roots, $\tilde{\phi}_1$ and $\tilde{\phi}_2$ \cite{Chopp2001}. 
	The smallest root which is larger than the surrounding stencil nodes is taken to be the accepted value. In the example case given the solution would be the smaller 
	of $\tilde{\phi}_1$ and $\tilde{\phi}_2$ that is larger than both $\phi_{i-1,j}$ and $\phi_{i,j+1}$. 
	
	The remaining grid nodes are updated in the following fashion. The grid node in the trial list with the smallest level set value is moved into the accepted list. 
	All nodes surrounding the newly accepted node which are in either the trial or distant lists are updated by solving Eq. (\ref{eq:StdFMM}). 
	Any updated nodes in the distant list
	are moved into the trial list. This procedure is repeated until no nodes remain in the trial list. The next node to be added to the accepted list is easily obtained if the trial list is 
	maintained as a sorted list such as a heap.

\section{The Augmented Fast Marching Method}
\label{sec:3.0}

	To improve the accuracy of fast marching based reinitialization schemes it is proposed to solve an extension of Eq. (\ref{eq:Eikonal}). 
	Begin by taking up to second order derivatives of the square of the Eikonal equation with $F=1$:
	\begin{align}
		\label{eq:afmm0}
		\nabla\phi \cdot\nabla\phi&=1, \\
		\label{eq:afmm1}
		\nabla\left(\nabla\phi \cdot\nabla\phi\right)&=0, \\
		\label{eq:afmm2}
		\nabla\nabla\left(\nabla\phi \cdot\nabla\phi\right)&=0.
	\end{align}
	In two dimensions this results in six equations, 
	\begin{align}
		\label{eq:2D_pde_phi} \phi_x^2+\phi_y^2&=1,\\
		\label{eq:2D_pde_phix} \phi_x \phi_{xx} +\phi_y \phi_{xy}&=0, \\
		\label{eq:2D_pde_phiy} \phi_x \phi_{xy} +\phi_y \phi_{yy}&=0, \\
		\label{eq:2D_pde_phixx} \phi_{xx}^2+\phi_{xy}^2+\phi_x \phi_{xxx}+\phi_y \phi_{xxy}&=0, \\
		\label{eq:2D_pde_phiyy} \phi_{yy}^2+\phi_{xy}^2+\phi_x \phi_{xyy}+\phi_y \phi_{yyy}&=0, \\
		\label{eq:2D_pde_phixy} \phi_{xx}\phi_{xy}+\phi_{yy}\phi_{xy}+\phi_x \phi_{xxy}+\phi_y \phi_{xyy}&=0,
	\end{align}
	with $\phi_x$ denoting partial derivative of $\phi$ with respect to $x$. The values of interest are the level set function and up to the second derivatives of the 
	level set function: $\phi$, $\phi_x$, $\phi_y$, $\phi_{xx}$, $\phi_{yy}$, and $\phi_{xy}$. By using finite difference approximations to first order derivatives
	a set of six equations is obtained for the six unknowns. The particular discretization used in this work is given in Sec. \ref{sec:3.2}, while the extension to three
	dimensions is presented in Sec. \ref{sec:3.3}.

\subsection{Initializing the AFMM}
\label{sec:3.1}

	Initialization of the augmented fast marching method proceeds in a manner similar to the standard fast marching method. Let the grid cell $\Omega_{i,j}$ contain
	the interface. It is possible to define a bicubic interpolant over this cell using the given data. As this work was designed to work with the gradient augmented
	level set method it is assumed that the level set value, $\phi$, and the gradient of the level set, $\nabla\phi=\vec{\psi}=\left(\psi^x,\psi^y\right)$,
	is available at the four grid 
	nodes associated with $\Omega_{i,j}$. To compute the bicubic interpolant it is necessary to define $\phi_{xy}$ at the four grid points. 
	This value is calculated as the average of the derivatives of the gradient field. At any point $\phi_{xy}=(D_x \psi^y+D_y \psi^x)/2$ where 
	$D_x$ and $D_y$ represent the centered second order finite difference approximations to the first derivative.
	
	To initialize the augmented fast marching method all nodes associated with grid cells containing the interface are moved into the accepted list. Values 
	for the level set, $\phi$, the gradient of the level set, $\vec{\psi}$, and the Hessian of the level set, $\vec{H}$, are calculated 
	at each of these initially accepted points.

\begin{figure}[!ht]
		\centering
		\includegraphics[width=0.4\textwidth]{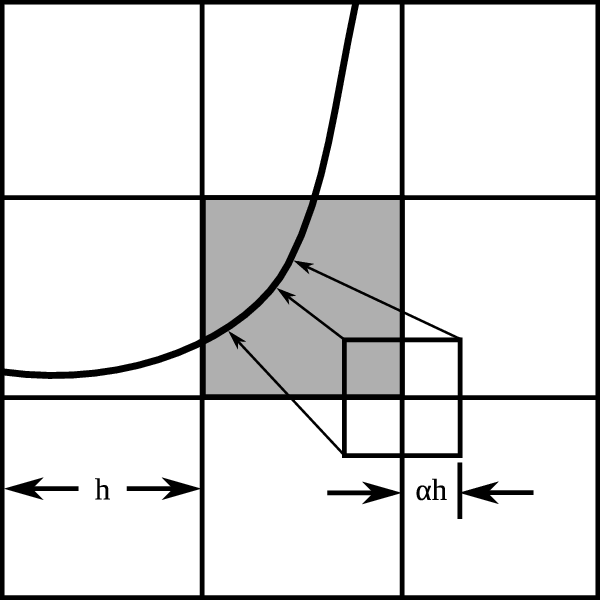}
		\label{fig:SubGrid}
		\caption{The initialization grid at a grid point $\vec{x}_{i,j}$ with a sub-grid centered at the node. 
		The distance from each node of the sub-grid to the interface is calculated and shown for three nodes.		
		The values of $\phi$, $\vec{\psi}$ and $\vec{H}$ are then computed using finite difference
		approximations on the sub-grid.}		
	\end{figure}

	The initialization follows a technique developed for accurately calculating the curvature of level set functions \cite{Macklin2006}.
	Let the grid node $\vec{x}_{i,j}$ be in the initially accepted list. A $3\times3$ sub-grid is centered at $\vec{x}_{i,j}$. The spacing of this sub grid is
	taken to be $\alpha h$, where $\alpha<1$ and $h$ the uniform grid 
	spacing, see Fig. \ref{fig:SubGrid}.
	Each of the nine points in the sub-grid have a signed distance function value calculated by minimizing the function $\|\vec{y}-\vec{x}_0\|^2_2$ 
	subject to $P(\vec{y})=0$, where	$P(\vec{x})$ is the bicubic interpolant over the grid cell $\Omega_{i,j}$ and $\vec{x}_0$ is the sub-grid point. 
	The required derivatives are then obtained by standard 
	second-order finite difference schemes using the sub-grid data. Due to the small sub-grid size the accuracy of the initialized values is extremely high, see the 
	results sections for details about the convergence rate.
	
\subsection{Updating Nodes in the AFMM}
\label{sec:3.2}

	The remaining grid nodes are updated in an ordered manner. Nodes adjacent to the initially accepted list are
	updated by solving the discretization of Eqs. (\ref{eq:afmm0})-(\ref{eq:afmm2}). These nodes are placed 
	into the trial list which is kept as a heap sort to ensure the fast retrial of the node with the smallest level set value. The node in the trial list with 
	the smallest value is placed into the accepted list and any non-accepted adjacent nodes are updated. 
	If an adjacent node is in the distant list it is moved into the trial list. This procedure is repeated until no nodes remain in the trial list.
	
	The updating procedure for the AFMM is more involved than the standard FMM. In the standard FMM a single nonlinear equation needs to be solved, Eq. (\ref{eq:StdFMM}).
	In the 2D AFMM there are six nonlinear equations. The solution obtained will depend on the particular discretization of the equations given in Sec. \ref{sec:3.0}. 
	In the case of two accepted nodes, one in each Cartesian direction, the discretization chosen is
	\begin{align}
		\label{eq:2D_phi} \psi^x (D_x^{\pm} \phi) + \psi^y (D_y^{\pm} \phi) &= 1, \\
		\label{eq:2D_phix} \psi^x (D_x^{\pm} \psi^x) + \psi^y (D_y^{\pm} \psi^x) & =0, \\
		\label{eq:2D_phiy} \psi^x (D_x^{\pm} \psi^y) + \psi^y (D_y^{\pm} \psi^y) & =0, \\
		\label{eq:2D_phixx} H^{xx} H^{xx} + H^{xy} H^{xy} + \psi^x (D_x^{\pm} H^{xx}) + \psi^y (D_y^{\pm} H^{xx}) &=0, \\
		\label{eq:2D_phiyy} H^{yy} H^{yy} + H^{xy} H^{xy} + \psi^x (D_x^{\pm} H^{yy}) + \psi^y (D_y^{\pm} H^{yy}) &=0, \\
		\label{eq:2D_phixy} H^{xx} H^{xy} + H^{yy} H^{xy} + \psi^x (D_x^{\pm} H^{xy}) + \psi^y (D_y^{\pm} H^{xy}) &=0.
	\end{align}	
	The operator $D_x^{\pm}$ and $D_y^{\pm}$ represent the appropriate one-sided derivatives at the grid point to be updated,
	see Sec. \ref{sec:2.3} or Ref. \cite{Chopp2001} . Also note that due to symmetry there are three components to the Hessian matrix which we denote as 
	$\phi_{xx}= H^{xx}$, $\phi_{yy}= H^{yy}$, and $\phi_{xy}=H^{xy}$.
	It is worth noting that this particular discretization allows for the solution of the level set and gradient field first, Eqs.
	(\ref{eq:2D_phi})-(\ref{eq:2D_phiy}). Once a valid solution set $\{\phi,\vec{\psi}\}$ is calculated 
	the solution of the Hessian field, Eqs. (\ref{eq:2D_phixx})-(\ref{eq:2D_phixy}) can be determined.
	
	To determine if a calculated solution to the gradient system, Eqs. (\ref{eq:2D_phi})-(\ref{eq:2D_phiy}), 
	is valid acceptance criteria similar to the classical FMM must be implemented. Due to the larger amount of information
	the acceptance criteria of the classical FMM is augmented. In addition to requiring that the solution level set value be larger than either of the grid point's
	neighbors it is also required that the calculated gradient be in the same general direction as the neighboring nodes. For example, let the nodes 
	$\vec{x}_{i+1,j}$ and $\vec{x}_{i,j+1}$ be accepted nodes when updating the value at $\vec{x}_{i,j}$. In this case the chosen finite difference approximations
	would be the stencils for the positive one-sided derivatives in each direction, $D_x^+$ and $D_y^+$. 
	Let $\tilde{\phi}$ and $\tilde{\vec{\psi}}=\left(\tilde{\psi}^x, \tilde{\psi}^y\right)$ be solutions to Eqs. (\ref{eq:2D_phi})-(\ref{eq:2D_phiy}). A valid solution set
	would satisfy $\tilde{\phi}\geq\phi_{i+1,j}$, $\tilde{\phi}\geq\phi_{i,j+1}$, $\tilde{\vec{\psi}}\cdot \vec{\psi}_{i+1,j}\geq 0$, and 
	$\tilde{\vec{\psi}}\cdot \vec{\psi}_{i,j+1}\geq 0$.
	If the solution set $\{\tilde{\phi},\tilde{\vec{\psi}}\}$ satisfies these conditions then it is taken as valid.

	When updating the values at $\vec{x}_{i,j}$ it is possible that only a single neighboring node is in the accepted list. In this case the discretization is 
	modified to account for the reduced amount of information. To use only information in the $x$-direction the system becomes
	\begin{align}
		\label{eq:1Dx_phi} \psi^x (D_x^{\pm} \phi) + \psi^y \psi^y &= 1, \\
		\label{eq:1Dx_phix} \psi^x (D_x^{\pm} \psi^x) + \psi^y (D_x^{\pm} \psi^y) & =0, \\
		\label{eq:1Dx_phiy} \psi^x (D_x^{\pm} \psi^y)  & =0, \\
		\label{eq:1Dx_phixx} H^{xx} H^{xx} + H^{xy} H^{xy} + \psi^x (D_x^{\pm} H^{xx}) + \psi^y (D_x^{\pm} H^{xy}) &=0, \\
		\label{eq:1Dx_phiyy} H^{yy} H^{yy} + H^{xy} H^{xy} + \psi^x (D_x^{\pm} H^{yy})  &=0, \\
		\label{eq:1Dx_phixy} H^{xx} H^{xy} + H^{yy} H^{xy} + \psi^x (D_x^{\pm} H^{xy}) + \psi^y (D_x^{\pm} H^{yy}) &=0.
	\end{align}		
	The general idea is to replace derivatives in the direction not present ($y$-direction in this case) with the equivalent $x$-direction derivatives.
	A similar formulation can be made for using information only from the $y$-direction and is given in Appendix \ref{sec:A}.

\subsection{Extension of the AFMM to 3D}
\label{sec:3.3}

	The extension to three dimensions is straightforward. Explicitly writing Eqs. (\ref{eq:afmm0})-(\ref{eq:afmm2}) results in ten equations,
	\begin{align}
		\label{eq:3D_pde_phi} \phi_x^2+\phi_y^2+\phi_z^2&=1,\\
		\label{eq:3D_pde_phix} \phi_x \phi_{xx} +\phi_y \phi_{xy} + \phi_z \phi_{xz}&=0, \\
		\label{eq:3D_pde_phiy} \phi_x \phi_{xy} +\phi_y \phi_{yy} + \phi_z \phi_{yz}&=0, \\
		\label{eq:3D_pde_phiz} \phi_x \phi_{xz} +\phi_y \phi_{yz} + \phi_z \phi_{zz}&=0, \\				
		\label{eq:3D_pde_phixx} \phi_{xx}^2+\phi_{xy}^2+\phi_{xz}^2+\phi_x \phi_{xxx}+\phi_y \phi_{xxy}+\phi_z \phi_{xxz}&=0, \\
		\label{eq:3D_pde_phiyy} \phi_{xy}^2+\phi_{yy}^2+\phi_{yz}^2+\phi_x \phi_{xyy}+\phi_y \phi_{yyy}+\phi_z \phi_{yyz}&=0, \\
		\label{eq:3D_pde_phizz} \phi_{xz}^2+\phi_{yz}^2+\phi_{zz}^2+\phi_x \phi_{xzz}+\phi_y \phi_{yzz}+\phi_z \phi_{zzz}&=0, \\				
		\label{eq:3D_pde_phixy} \phi_{xx}\phi_{xy}+\phi_{yy}\phi_{xy}+\phi_{xz}\phi_{yz}+\phi_x \phi_{xxy}+\phi_y \phi_{xyy}+\phi_z \phi_{xyz}&=0, \\
		\label{eq:3D_pde_phixz} \phi_{xx}\phi_{xz}+\phi_{zz}\phi_{xz}+\phi_{xy}\phi_{yz}+\phi_x \phi_{xxz}+\phi_y \phi_{xyz}+\phi_z \phi_{xzz}&=0, \\
		\label{eq:3D_pde_phiyz} \phi_{yy}\phi_{yz}+\phi_{zz}\phi_{yz}+\phi_{xz}\phi_{xy}+\phi_x \phi_{xyz}+\phi_y \phi_{yyz}+\phi_z \phi_{yzz}&=0.	
	\end{align}	
	The quantities of interest are the level set function $\phi$ and derivatives of $\phi$ up to second order. Using upwind finite difference approximations for the 
	first derivative this results in a set of ten nonlinear equations for the ten unknowns.		
	
	In the three dimensional case initialization of a grid point $\vec{x}_{i,j,k}$ occurs over a $3\times 3\times 3$ sub-grid centered on $\vec{x}_{i,j,k}$. To use numerical approximations of the
	derivatives it is not necessary to determine the level set value for all 27 points on this sub-grid as only 19 points are used 
	during the derivative calculations. As in the two-dimensional case these sub-grid points have their distance to the interface calculated by minimizing 
	$\|\vec{y}-\vec{x}_0\|^2_2$ subject to $P(\vec{y})=0$, where
	$P(\vec{x})$ is the tricubic interpolant over the grid cell $\Omega_{i,j,k}$. When computing the tricubic coefficients it is sufficient to ensure that
	the following values are satisfied at the eight corner nodes of the grid cell: $\phi$, $\phi_x$, $\phi_y$, $\phi_z$, $\phi_{xy}$, $\phi_{xz}$, $\phi_{yz}$, and
	$\phi_{xyz}$.	
	Assuming that only the level set and the gradient of the level set are
	available during reinitialization any higher order derivatives are obtained by averaging the appropriate derivatives. For example 
	the third-order derivative necessary for the tricubic function is given by $\phi_{xyz}=(D_{yz} \psi^x+D_{xz} \psi^y+D_{xy} \psi^z)/3$, where $D_{xy}$ is the 
	finite difference approximation to $\partial_{xy}$.
	
	Updating the remaining grid points proceeds in the same fashion as the two-dimensional case.
	Due to the higher number of dimensions there are seven possibilities from where information will be propagating. The specific discretizations for
	all possibilities have been presented in Appendix \ref{sec:B}.

\subsection{The AFMM Algorithm}
\label{sec:3.4}
	One consequence of choosing the above discretizations is that the solution for the gradient system, Eqs. (\ref{eq:2D_phi})-(\ref{eq:2D_phiy}) for the two dimensional case,
	can be computed first. Once a valid solution set $\{\phi,\vec{\psi}\}$ is calculated 
	the solution for the Hessian field, Eqs. (\ref{eq:2D_phixx})-(\ref{eq:2D_phixy}) can be determined. In practice this is done by first computing for the gradient
	field in the entire domain. This not only results in level set and gradient functions at every grid point, it also determines the proper ordering to ensure
	that upwind derivatives are calculated correctly. The Hessian field is then obtained by solving Eqs. (\ref{eq:2D_phixx})-(\ref{eq:2D_phixy}) using this pre-determined 
	ordering. 
	
	The algorithm for the augmented fast marching method can be summarized thusly:
	
	\begin{enumerate}
		\item{Mark all nodes as in the Distant list.}
		\item{Initialize all nodes associated with cells containing the interface by explicitly solving for the signed distance function on a sub-grid centered at the node. Derivatives 
			of the level set are calculated by standard finite difference approximations on this subgrid. Move these nodes into the Accepted list.}
		\item{For all nodes in the Distant list that lie next to a node in the Accepted list calculated updated values by solving Eqs.
		(\ref{eq:2D_pde_phi})-(\ref{eq:2D_pde_phiy}) in 2D or Eqs. (\ref{eq:3D_pde_phi})-(\ref{eq:3D_pde_phiz}) in 3D. Move these nodes into the Trial list.}
		\item{Select the node with the smallest level set value in the Trial list and move it into the Accepted list. All nodes next to the newly accepted grid point in either the Trial or Distant list
		have updated values calculated by solving Eqs. (\ref{eq:2D_pde_phi})-(\ref{eq:2D_pde_phiy}) in 2D or Eqs. (\ref{eq:3D_pde_phi})-(\ref{eq:3D_pde_phiz}) in 3D.
		Any node that is updated and in the Distant list is moved into the Trial list.}
		\item{Repeat Step 4 until the Trial list is empty.}
		\item{Calculate the Hessian field by solving Eqs. (\ref{eq:2D_pde_phixx})-(\ref{eq:2D_pde_phixy}) in 2D 
			or Eqs. (\ref{eq:3D_pde_phixx})-(\ref{eq:3D_pde_phiyz}) in 3D. Update nodes using the same order as they were added to the Accepted list.}		
	\end{enumerate}

\section{Two Dimensional Results}
\label{sec:4.0}

	Here convergence and sample results are presented for two dimensional interfaces. All numerical derivatives were calculated using standard first-order one-sided finite 
	difference schemes.
	The domain is the region $[-2,2]^2$ and a uniform grid 
	spacing $h$ is used. No domain boundary conditions are needed as all information flows from the interface outwards.	
	In all cases the the interface is initially described by a level set function $\phi_0$ and its gradient field, $\vec{\psi}_0=\nabla\phi_0$.
	The sub-grid is taken to have a spacing of $0.1h$. The SLSQP algorithm  \cite{Kraft1988,Kraft1994} of the NLopt software library \cite{Johnson} was 
	used to determine the closest point on the interface during initialization. All nonlinear systems were solved using the GNU Scientific Library \cite{Galassi2003}.

\subsection{Investigation of the Expected Errors for the Gradient System}
\label{sec:expected_errors}
	
\begin{figure}[!ht]
		\centering
		\includegraphics[width=0.4\textwidth]{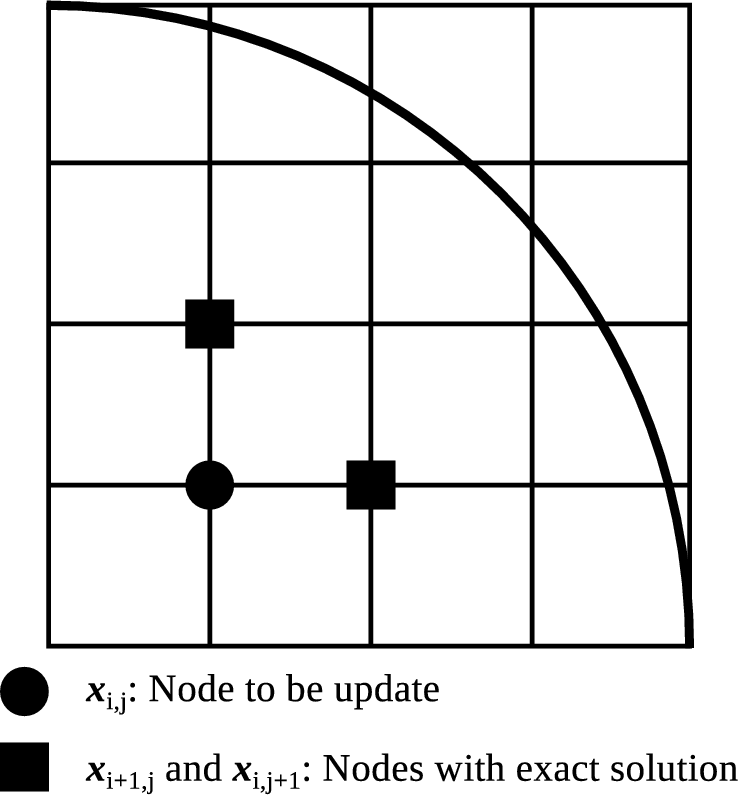}
		\caption{Sample stencil used for analysis of expected errors. The grid point $\vec{x}_{i,j}$ is to be updated by the AFMM using the 
		nodes $\vec{x}_{i+1,j}$ and $\vec{x}_{i,j+1}$. The nodes $\vec{x}_{i+1,j}$ and $\vec{x}_{i,j+1}$
		have the exact solution prescribed.}
		\label{fig:AnalyticAFMM}
	\end{figure}

	The accuracy of a fast marching method depends on the order in which nodes are updated. 
	A standard error analysis is difficult due to the nonlinear nature of the systems involved. Instead sample analytic solutions for a circular interface of radius $R_0$ are considered.
	This interface is to be reinitialized on a grid with a uniform spacing of $h$. Consider the AFMM updating of 
	a grid point given by $\vec{x}_{i,j}=\left(x,x\right)$, $x>0$, and $\sqrt{2x^2}<R_0$. This point lies on the domain diagonal.
	In this situation it is known that the grid points to the right, $\vec{x}_{i+1,j}$, and above, $\vec{x}_{i,j+1}$, are closer to the interface than $\vec{x}_{i,j}$, Fig.
	\ref{fig:AnalyticAFMM}. If the neighboring nodes have the exact solution prescribed the error of node $\vec{x}_{i,j}$ can be investigated.
	
	Consider the gradient system, Eqs. (\ref{eq:2D_pde_phi})-(\ref{eq:2D_pde_phiy}). Using the exact solution for nodes $\vec{x}_{i+1,j}$ and $\vec{x}_{i,j+1}$
	the level set, $\phi$, and gradient vector, $\nabla\phi=\vec{\psi}=\left(\psi^x,\psi^y\right)$, at $\vec{x}_{i,j}$ can be calculated as 
	\begin{align}
		\label{eq:exactphi} \phi=&\frac{2 x\sqrt{h^2+2hx+2x^2}}{h+2x}-R_0, \\
		\label{eq:exactpsi} \psi^x=\psi^y=&\frac{h+2x}{2\sqrt{h^2+2hx+2x^2}}.
	\end{align}
		
\begin{figure}[!ht]
		\centering
		\subfigure[Level Set Error]{\includegraphics[width=0.475\textwidth]{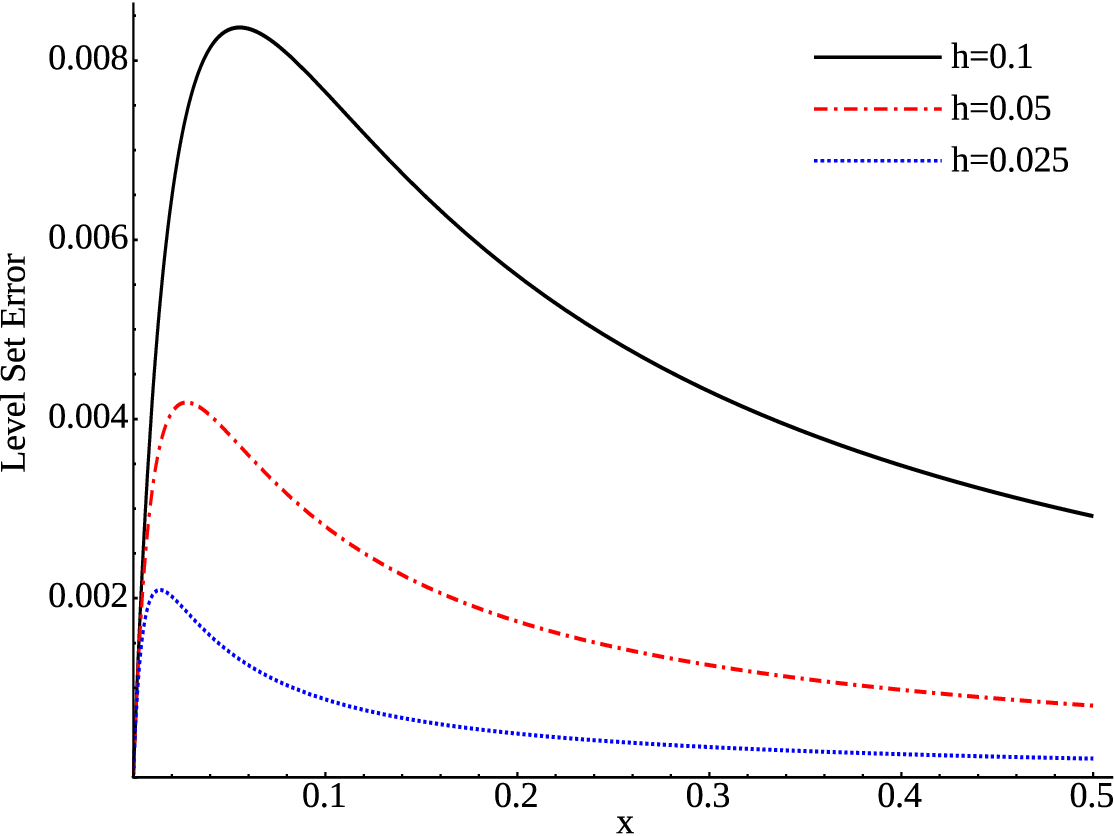} \label{fig:AFMM_Error_LevelSet}} 
		\subfigure[Gradient Error]{\includegraphics[width=0.475\textwidth]{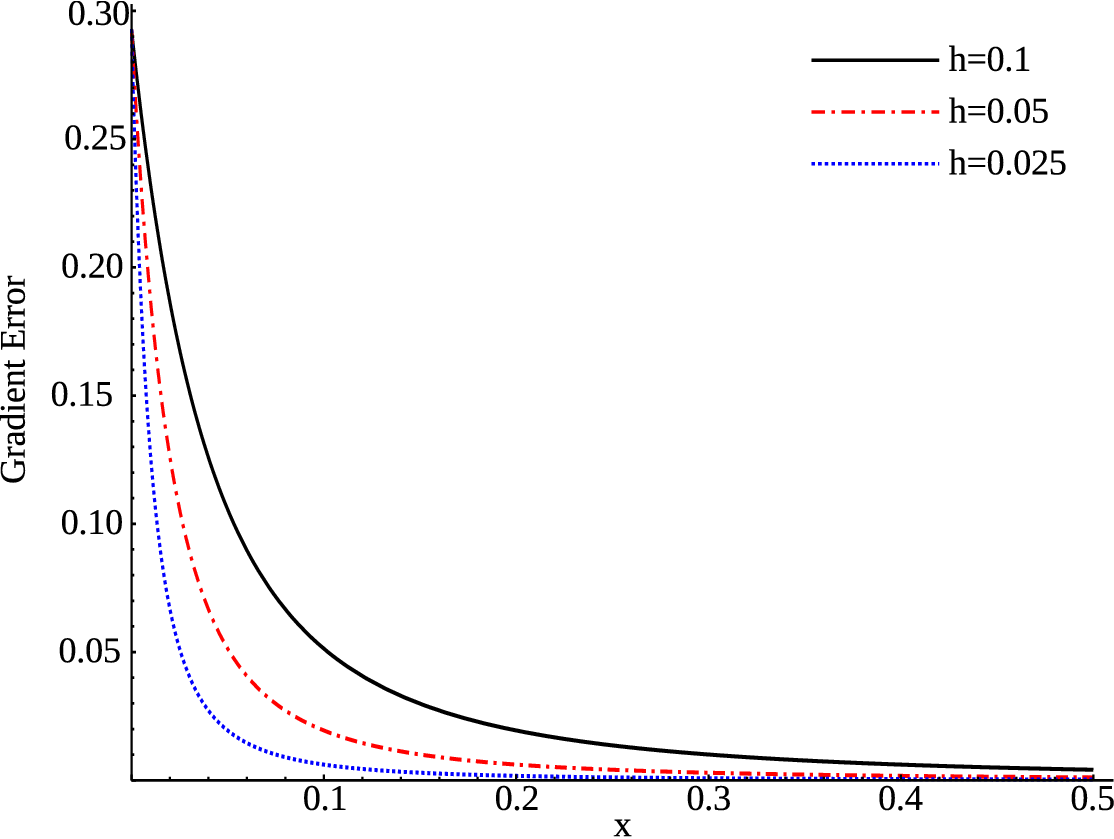} \label{fig:AFMM_Error_Gradient}}
		\caption{The error of the level set and gradient functions based on the stencil given in Fig. \ref{fig:AnalyticAFMM} for grid spacings of $h=0.1,0.05,\textrm{and }0.025$.}
		\label{fig:AFMM_Error}
	\end{figure}

	The true solutions are $\phi_{true}=\sqrt{2x^2}-R_0$ and $\psi_{true}^x=\psi_{true}^y=1/\sqrt{2}$. The error of the level set is calculated as 
	$\phi-\sqrt{2x^2}+R_0$ while the error of the gradient vector is $\sqrt{\left(\psi^x-1/\sqrt{2}\right)^2+\left(\psi^y-1/\sqrt{2}\right)^2}$, where 
	$\phi$, $\psi^x$, and $\psi^y$ are given in Eqs. (\ref{eq:exactphi}) and (\ref{eq:exactpsi}). These errors are shown in Fig. \ref{fig:AFMM_Error} for 
	locations ranging from $x=0$ to $x=0.5$ and grid spacings of $h=0.1, 0.05, \textrm{and }0.025$.

	From  the results shown in Fig. \ref{fig:AFMM_Error} it becomes apparent that as the grid size decreases the overall error decreases for both the level set and gradient
	vector. The maximum error of the level set function observes a first-order convergence. The gradient vector, on the other hand, has a fixed error of 
	$\left(\sqrt{2}-1\right)/\sqrt{2}$ as $x$ approaches zero irregardless of the grid spacing. This result should not be unexpected. As $x$ approaches the origin
	the variation of the gradient field increases inversely to the distance from the origin. For example, given a circular interface and a point $\vec{x}=(x,x)$ the exact variation of 
	the x component of the gradient field in the x-direction is given by $\partial_x\psi^x = 1/\left(2\sqrt{2}x\right)$.	
	As the distance from the origin decreases due to a
	smaller grid spacing the variation in the gradient field will increase by an inverse amount. This results in a fixed error being introduced into the system near the origin.
	Despite this $\mathcal{O}(1)$ error of the gradient field the overall error decreases rapidly. 
	As the Hessian system depends on the solution of the gradient vector it should be expected that the errors for the Hessian field, and thus the curvature,
	will increase as one approaches the origin.
	
\begin{figure}[!ht]
		\centering
		\subfigure[Level Set Error]{\includegraphics[width=0.475\textwidth]{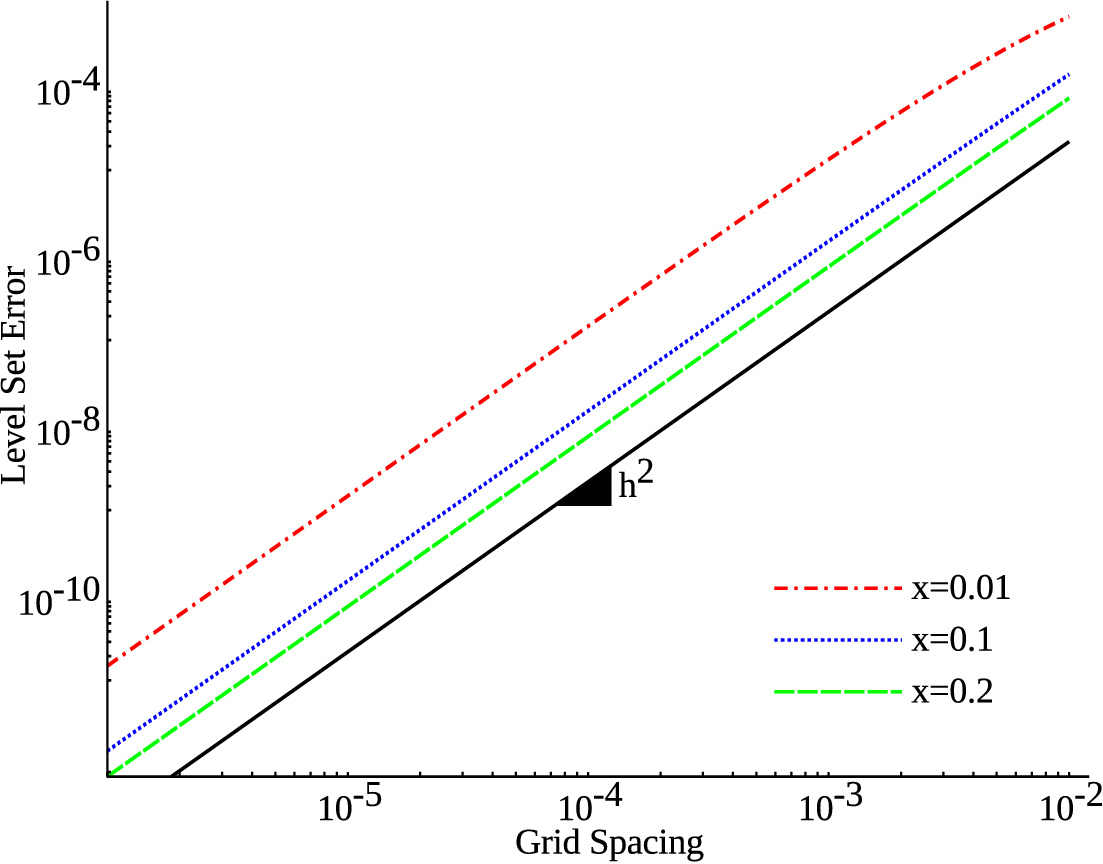} \label{fig:AFMM_LevelSet_Point_Error}}
		\subfigure[Gradient Error]{\includegraphics[width=0.475\textwidth]{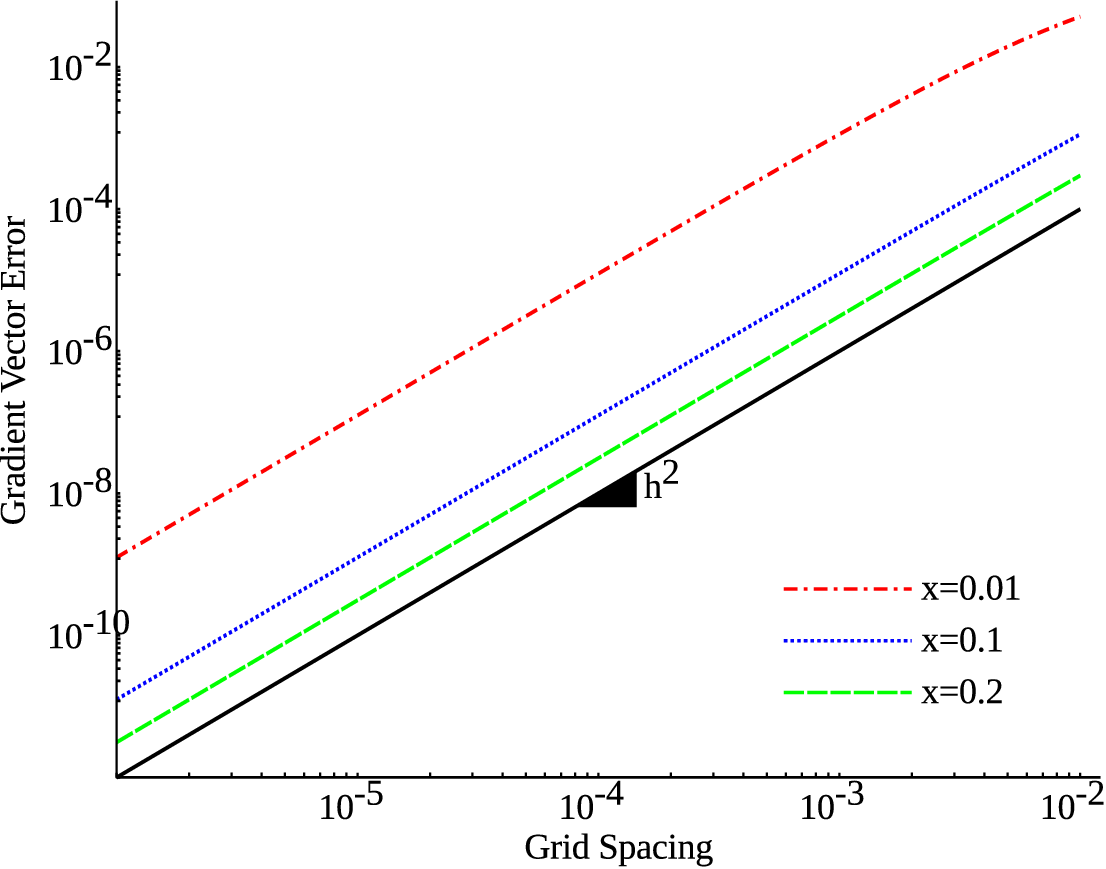} \label{fig:AFMM_Gradient_Point_Error}}
		\caption{The error of the level set and gradient functions based on the stencil given in Fig. \ref{fig:AnalyticAFMM} at the points $x=0.01, 0.1, \textrm{and }0.2$ for
		grid spacings ranging from $10^{-6}$ to $10^{-2}$.}		
		\label{fig:AFMM_point_error}
	\end{figure}

	To conclude this section it should be noted that the error at any given point will decrease for both the level set and gradient vector. As an example the error for 
	the level set and gradient vector are shown in Fig. \ref{fig:AFMM_point_error} for the points $x=0.01, 0.1, \textrm{and }0.2$ using grid spacing ranging from
	$10^{-6}$ to $10^{-2}$. For the three points considered second-order convergence is observed for both the level set and gradient vector solutions.
	
\subsection{Accuracy of the Initialization Method}
\label{sec:4.1}

\begin{figure}[!ht]
		\begin{center}
			\subfigure[$L_2$ Norm]{
				\includegraphics[width=0.475\textwidth]{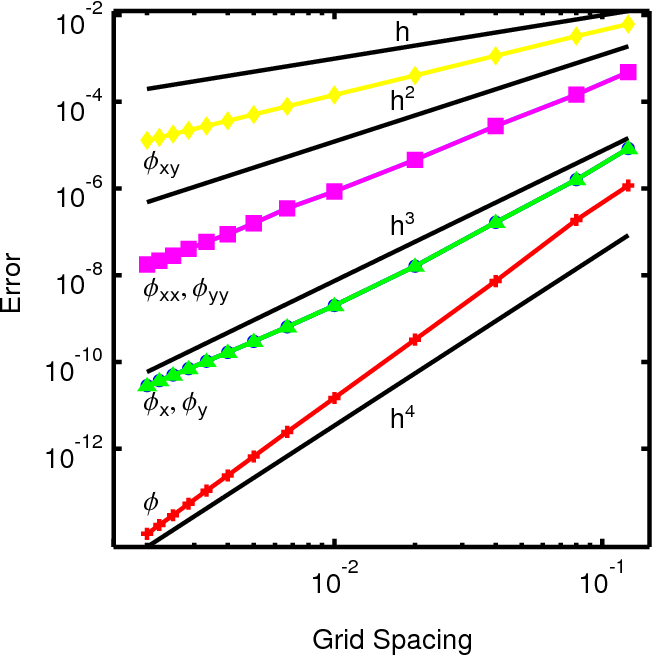}
				\label{fig:circle_initial_l2}
			}
			\subfigure[$L_\infty$ Norm]{
				\includegraphics[width=0.475\textwidth]{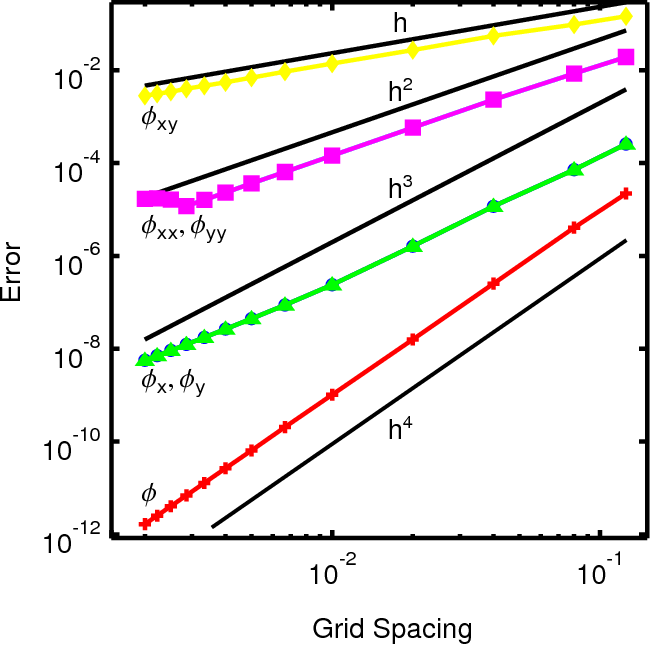}
				\label{fig:circle_initial_li}
			}
		\end{center}
		\caption{Convergence of the initialization procedure for a unit circle interface. The solid lines without symbols are convergence rates. }
		\label{fig:circle_initial}
	\end{figure}

	The accuracy of the initialization scheme presented in Sec \ref{sec:3.1} is shown here. 
	A circle of radius 1 with an initial level set given by 
	$\phi_0=e^{x^2+y^2}-e$ is considered. The resulting convergence rate is seen in Fig. \ref{fig:circle_initial}

	It is observed that the convergence rate for the level set is approximately $4^{th}$-order, for the gradient field it is $3^{rd}$ order, for $\phi_{xx}$ and $\phi_{yy}$ 
	$2^{nd}$-order convergence is seen, and the $\phi_{xy}$ value is $1^{st}$-order accurate.

\subsection{Accuracy of the 2D AFMM}
\label{sec:4.2}

\begin{figure}[!ht]
		\begin{center}
			\subfigure[$20^2$]{
				\includegraphics[width=0.475\textwidth]{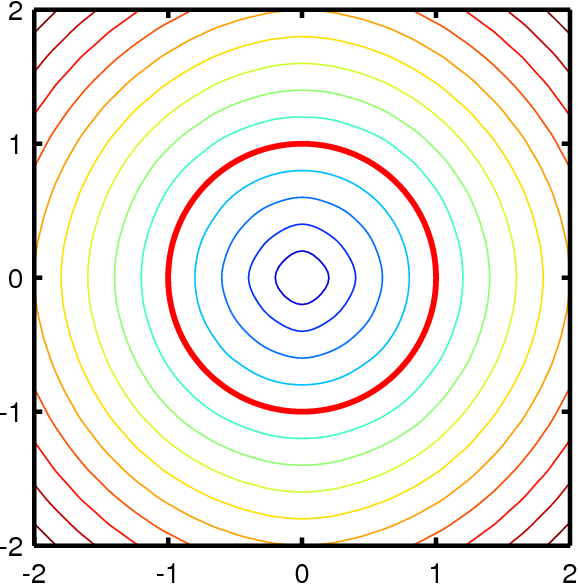}
			}
			\subfigure[$50^2$]{
				\includegraphics[width=0.475\textwidth]{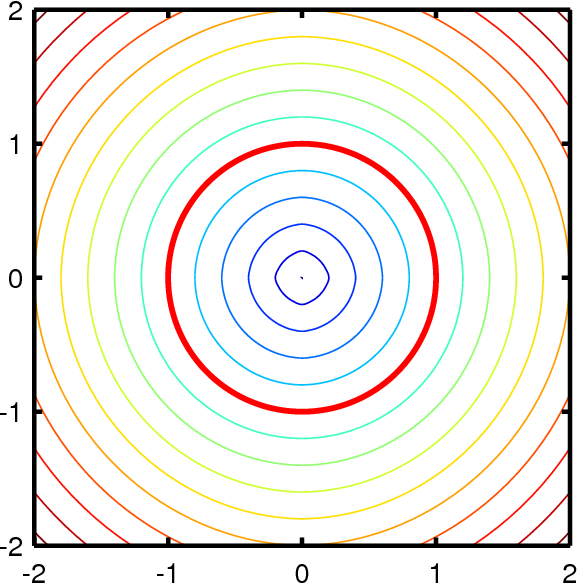}
			}\\
			\subfigure[$200^2$]{
				\includegraphics[width=0.475\textwidth]{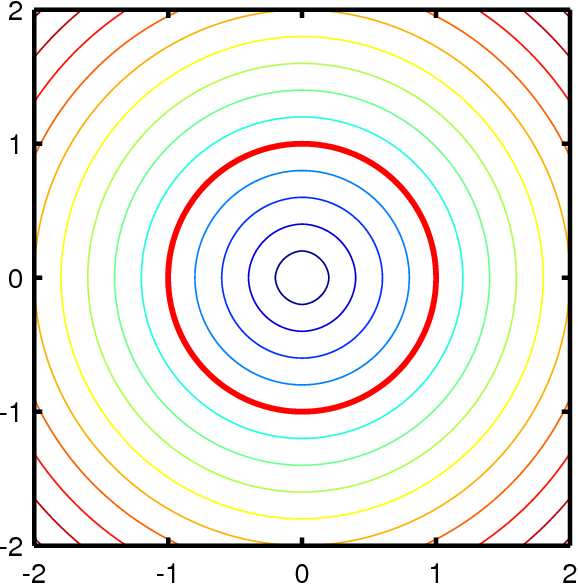}
			}
			\subfigure[$500^2$]{
				\includegraphics[width=0.475\textwidth]{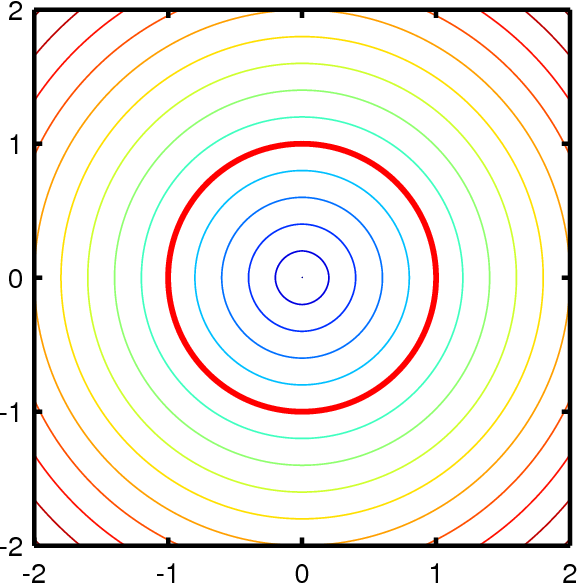}
			}\\
		\end{center}
		\caption{The level set after reinitialization for a unit circle with an initial level set of $\phi_0=e^{x^2+y^2}-e$. The interface is given by the thick red contour.
		These contours utilize the sub-grid information provided by the additional derivative information.}
		\label{fig:circle_phi}
	\end{figure}

\begin{figure}[!ht]
		\begin{center}
			\subfigure[$20^2$]{
				\includegraphics[width=0.475\textwidth]{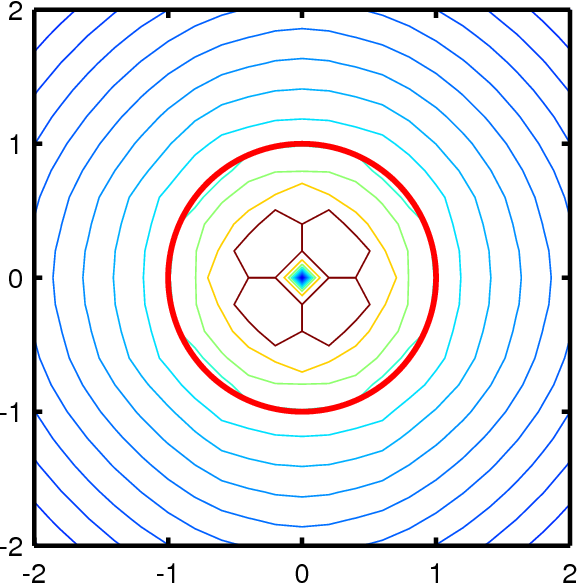}
			}
			\subfigure[$50^2$]{
				\includegraphics[width=0.475\textwidth]{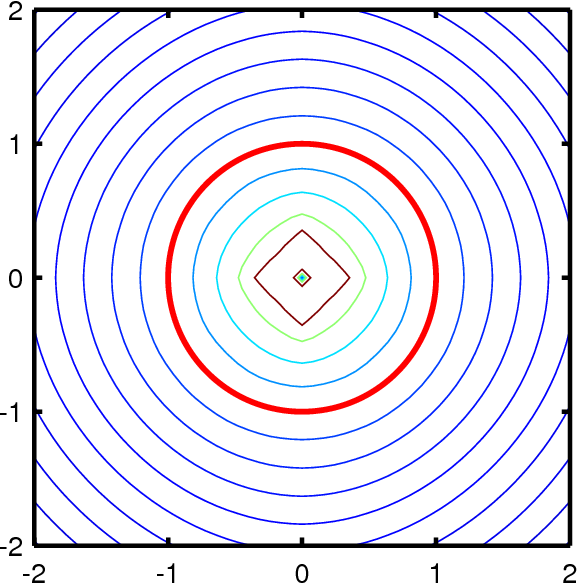}
			}\\
			\subfigure[$200^2$]{
				\includegraphics[width=0.475\textwidth]{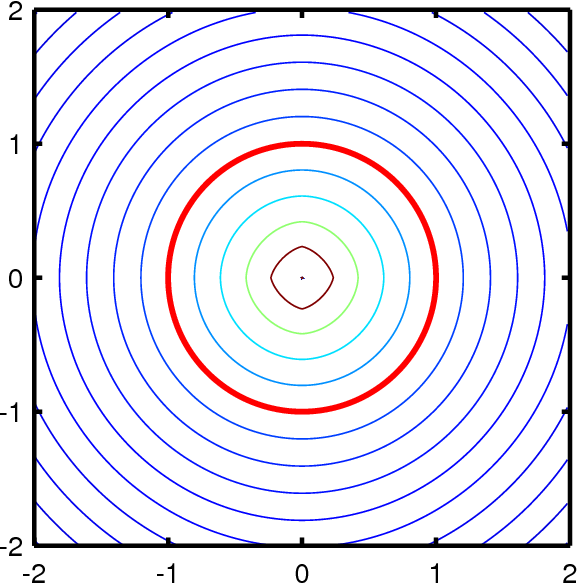}
			}
			\subfigure[$500^2$]{
				\includegraphics[width=0.475\textwidth]{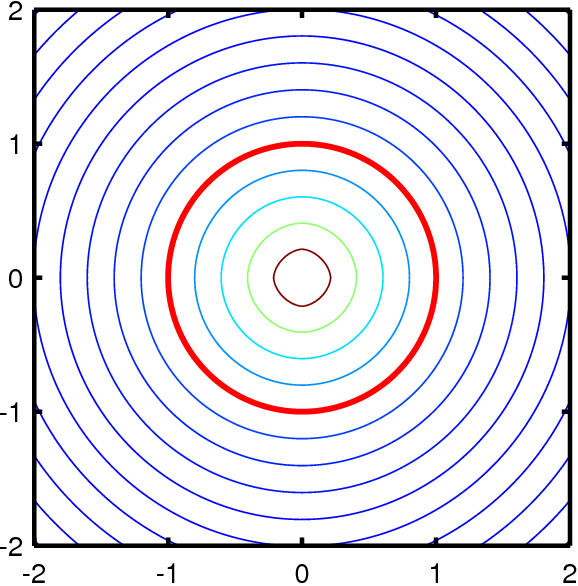}
			}\\
		\end{center}
		\caption{The curvature field after reinitialization for a unit circle with an initial level set of $\phi_0=e^{x^2+y^2}-e$. The curvature of 
		the interface is given by the thick red contour.}	
		\label{fig:circle_k}
	\end{figure}

	First consider a unit circle with the initial level set of $\phi_0=e^{x^2+y^2}-e$. Sample level set and curvature results at grid sizes ranging from $20^2$ to $500^2$ are shown in
	Figs. \ref{fig:circle_phi} and \ref{fig:circle_k}. The level set contours shown utilize the additional information provided by knowledge of the derivatives of the level 
	set function. Even in the coarsest mesh, $20^2$, the level set is extremely smooth and the curvature field is smooth outside of the interface. At such a coarse mesh
	there are not enough grid points to accuractly describe the large variations of curvature which occurs in the region given by $\phi<0$. As the number of grid points increases this error vanishes.	
	
\begin{figure}[!ht]
		\begin{center}
			\subfigure[$L_2$ Error]{
				\includegraphics[width=0.475\textwidth]{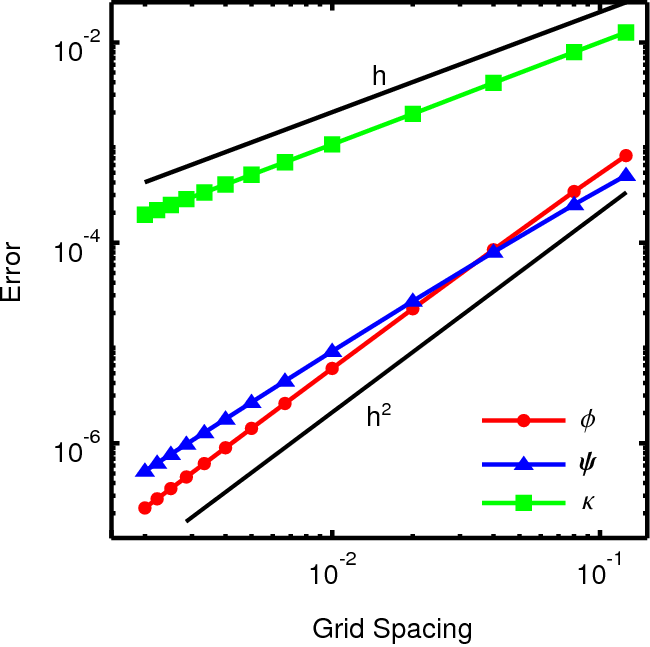}
			}
			\subfigure[$L_\infty$ Error]{
				\includegraphics[width=0.475\textwidth]{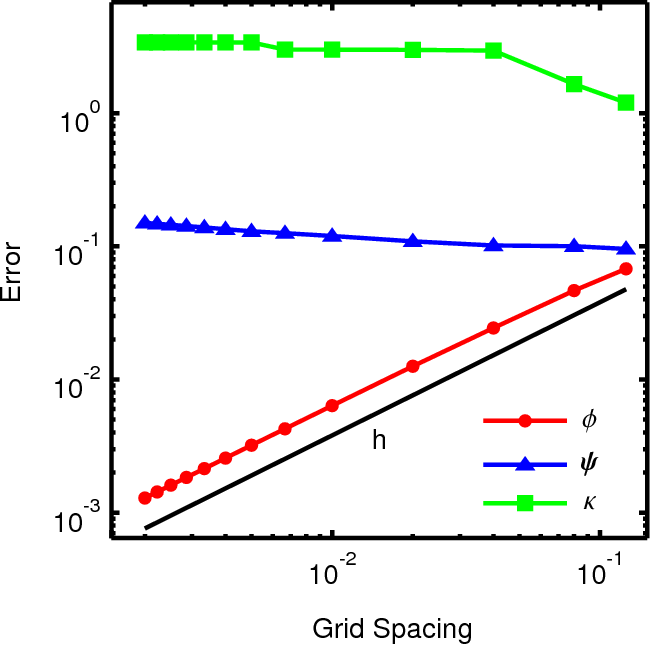}
			}
		\end{center}
		\caption{The error in the entire domain for the unit circle. The $L_\infty$ error does not converge for the gradient and curvature
		fields, as expected by the results of Sec. \ref{sec:expected_errors}.}
 		\label{fig:circle_entire_convergence}
	\end{figure}

\begin{figure}[!ht]
		\begin{center}
			\subfigure[$L_2$ Error]{
				\includegraphics[width=0.475\textwidth]{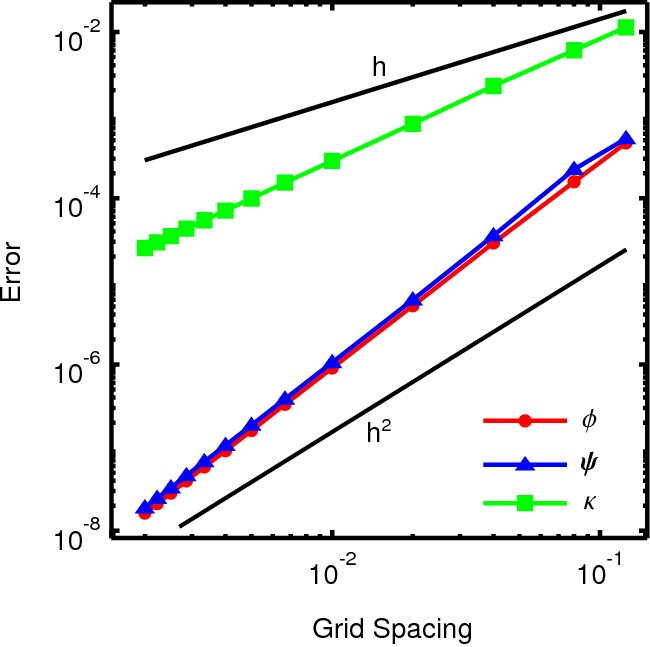}
			}
			\subfigure[$L_\infty$ Error]{
				\includegraphics[width=0.475\textwidth]{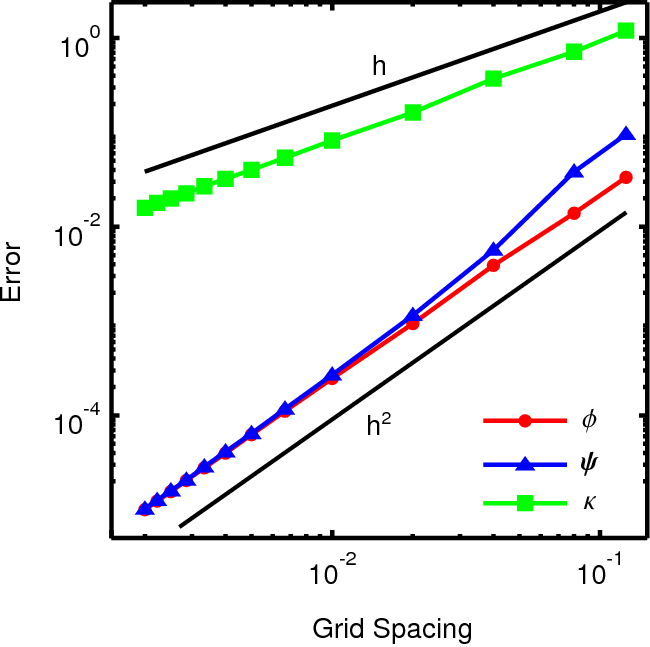}
			}
		\end{center}
		\caption{The error for the unit circle in the region given by $|\phi|\leq 9h$, where $h$ is the grid spacing. In this case
			both the $L_2$ and $L_\infty$ error-norms converge.}		
		\label{fig:circle_tube_convergence}
	\end{figure}

	Considering the results shown in Sec. \ref{sec:expected_errors} errors will be reported for both the entire domain and for the local region around the interface given by
	$|\phi|\leq 9h$. This width was chosen based on the stencils needed for a $5^{th}$-order WENO local level set scheme \cite{peng1999pde}. 
	To account for the fact that the curvature grows as $1/\sqrt{x^2+y^2}$ as one approaches the center of the given level set function all 
	curvature errors are reported as $(\kappa-\kappa_e)/\kappa_e$ where $\kappa$ is the curvature using the AFMM and $\kappa_e$ is the exact curvature. 
	The errors for the unit circle are given in Figs. \ref{fig:circle_entire_convergence} and \ref{fig:circle_tube_convergence}.
	
	Despite the use of first-order finite difference approximations the 
	resulting level set is second-order accurate and the gradient is $\mathcal{O}(h^{3/2})$ in the entire domain while the 
	curvature is $\mathcal{O}(h)$. When only considering the region close to the interface this increases
	to $\mathcal{O}(h^{5/2})$ for both the level set and gradient fields and $\mathcal{O}(h^{3/2})$ for the curvature field.
	
	The convergence results verify the conclusion obtained in Sec. \ref{sec:expected_errors}. The large variation in the gradient vector as one approaches the center of the circle 
	creates a uniform error that can not be overcome. The fact that the $L_2$-norm does converge indicates that this error is localized around the origin. This is further supported 
	by the convergence of both the $L_2$ and $L_\infty$-norms in the region given by $|\phi|\leq 9h$. 
	
\begin{figure}[!ht]
		\begin{center}
			\subfigure[$20^2$]{
				\includegraphics[width=0.475\textwidth]{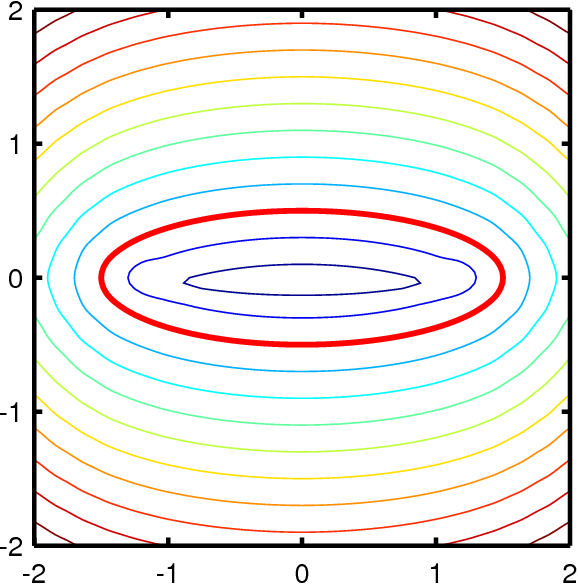}
			}
			\subfigure[$50^2$]{
				\includegraphics[width=0.475\textwidth]{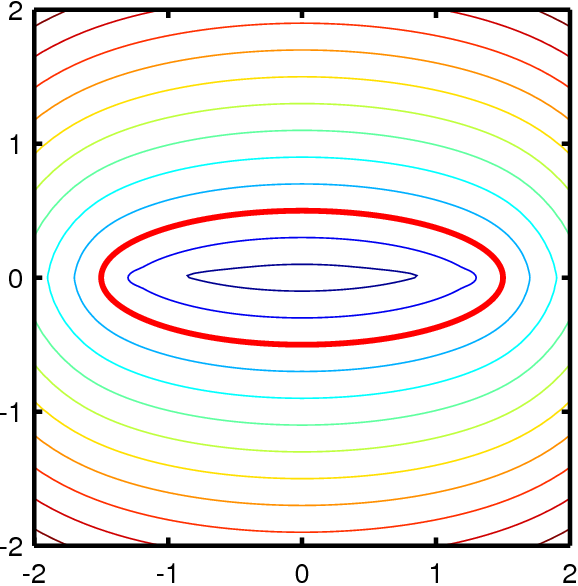}
			}\\
			\subfigure[$200^2$]{
				\includegraphics[width=0.475\textwidth]{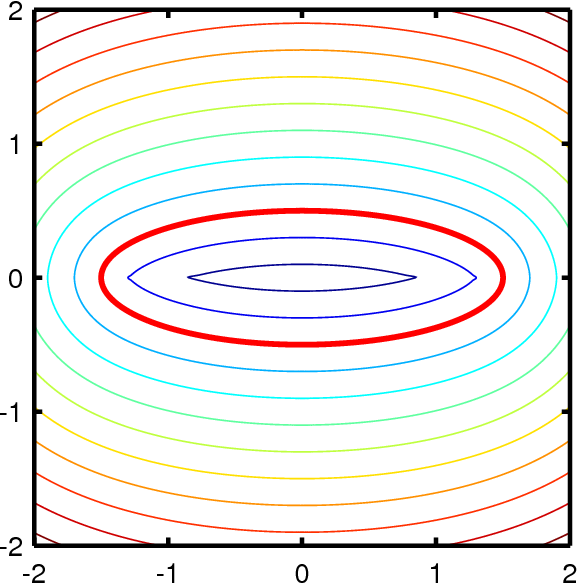}
			}
			\subfigure[$500^2$]{
				\includegraphics[width=0.475\textwidth]{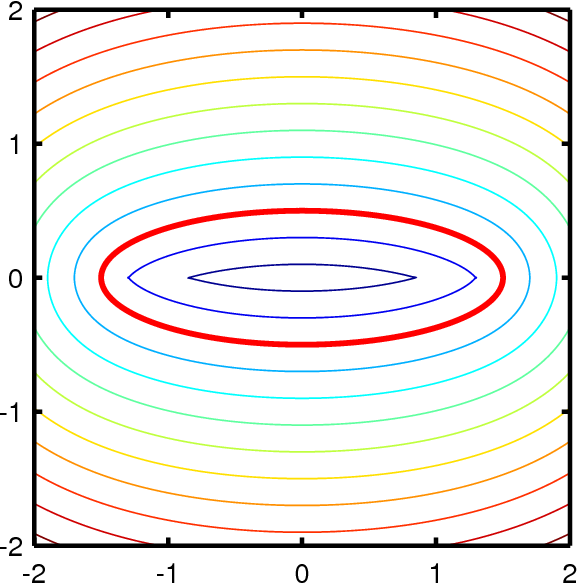}
			}\\
		\end{center}
		\caption{The level set after reinitialization for an elliptical interface with an initial level set of $\phi_0=(x/1.5)^2+(y/0.5)^2-1$. 
		The interface is given by the thick red contour.
		These contours utilize the sub-grid information provided by the additional derivative information.}
		\label{fig:ellipse_phi}
	\end{figure}

\begin{figure}[!ht]
		\begin{center}
			\subfigure[$20^2$]{
				\includegraphics[width=0.475\textwidth]{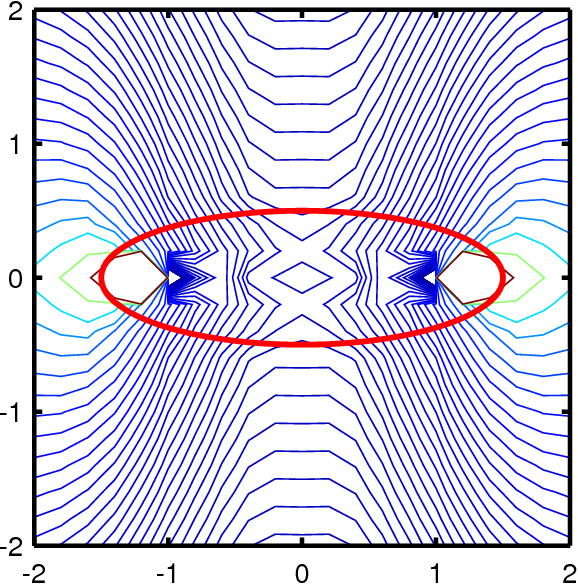}
			}
			\subfigure[$50^2$]{
				\includegraphics[width=0.475\textwidth]{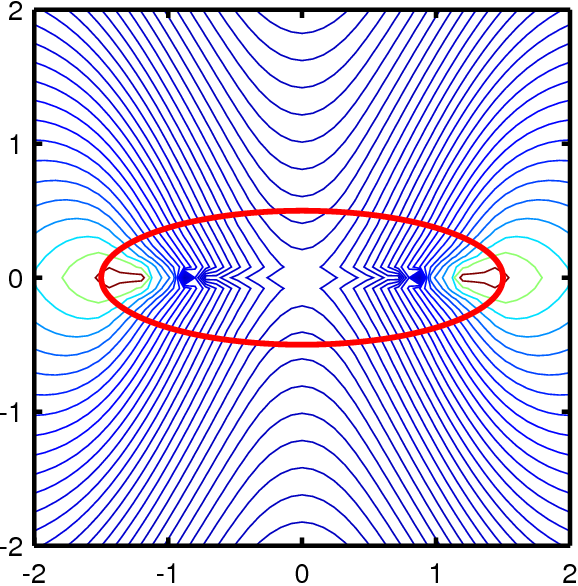}
			}\\
			\subfigure[$200^2$]{
				\includegraphics[width=0.475\textwidth]{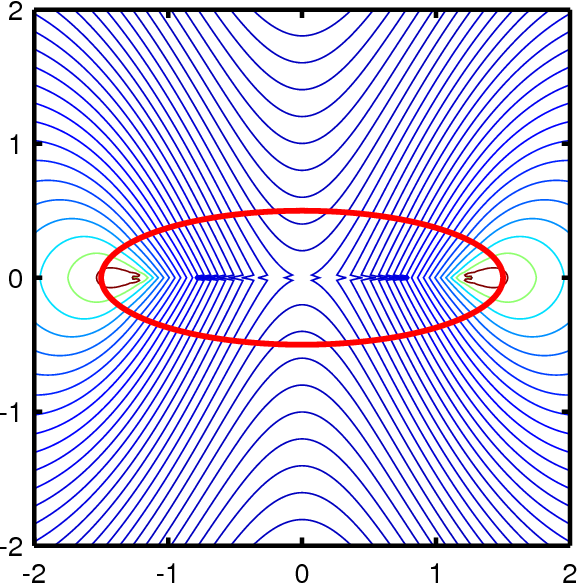}
			}
			\subfigure[$500^2$]{
				\includegraphics[width=0.475\textwidth]{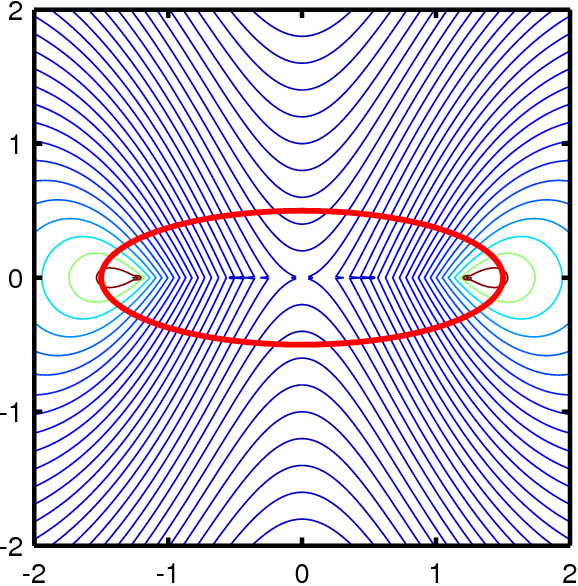}
			}\\
		\end{center}
		\caption{The curvature field after reinitialization for an elliptical interface with an initial level set of $\phi_0=(x/1.5)^2+(y/0.5)^2-1$.
		The curvature of the interface is given by the thick red contour.}	
		\label{fig:ellipse_k}
	\end{figure}

\begin{figure}[!ht]
		\begin{center}
			
			\includegraphics[width=0.475\textwidth]{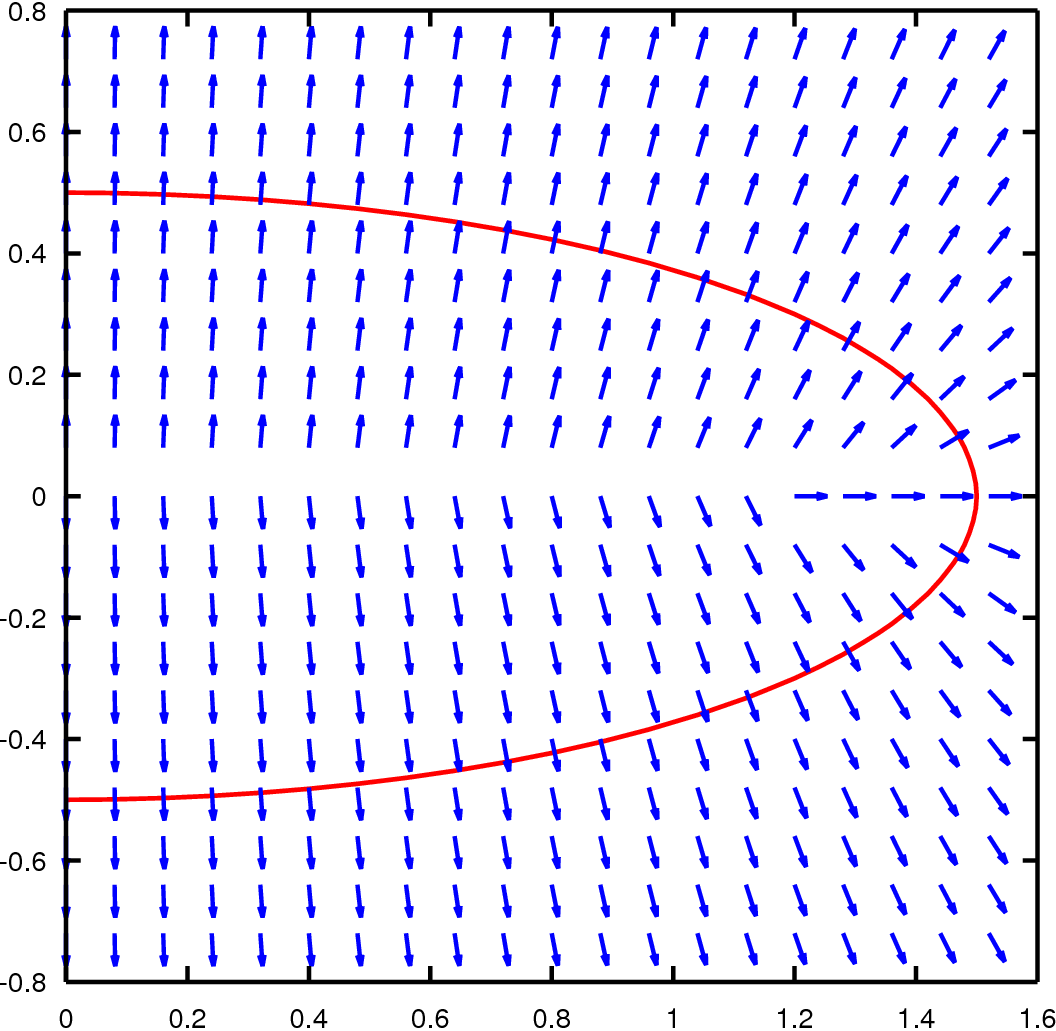}
			
		\end{center}
		\caption{Vector field of the reinitialized elliptical interface solution on a $50^2$ grid. A sharp discontinuity of the vector field occurs along the x-axis.}		
		\label{fig:ellipse_vectors}
	\end{figure}

	Next consider an elliptical interface with a major axis of 1.5 and a minor axis of 0.5: $\phi_0=(x/1.5)^2+(y/0.5)^2-1$. Results for various grid sizes are presented
	in Figs. \ref{fig:ellipse_phi} and \ref{fig:ellipse_k}. As for the unit circle the level set field is smooth even at coarse meshes. As the grid becomes finer the 
	level set field becomes smoother, particularly along the long axis of the ellipse. The curvature field is also smooth, even for the extremely coarse grid. For all of the grids
	the curvature field has the same qualitative shape. An issue with the curvature can be observed towards the center of the ellipse. This is due to the gradient vector field having
	a sharp discontinuity in that region. This is illustrated in Fig. \ref{fig:ellipse_vectors}, which shows the right half of the reinitialized solution on a $50^2$ grid.
	The resulting gradient vector field switches sign when crossing the x-axis, resulting in errors being introduced into the Hessian, and therefore the curvature, solution.

	As in the unit circle case the convergence of the elliptical interface is presented for both the entire domain and the region given by $|\phi|\leq 9h$. The exact
	solution was calculated by explicitly determining the signed distance function on a $4000^2$ grid.	
	The convergence 
	results	for the elliptical case match those of the unit circle.
	
	\begin{figure}[!ht]
		\begin{center}
			\subfigure[$L_2$ Error]{
				\includegraphics[width=0.475\textwidth]{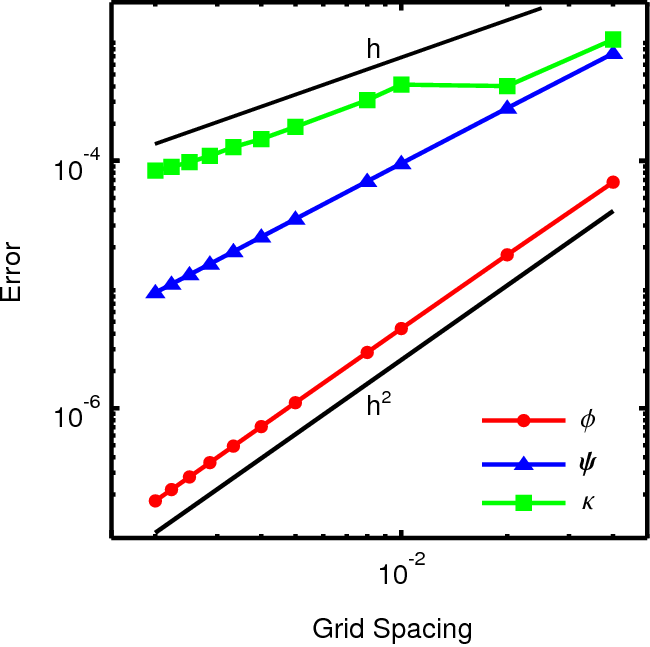}
			}
			\subfigure[$L_\infty$ Error]{
				\includegraphics[width=0.475\textwidth]{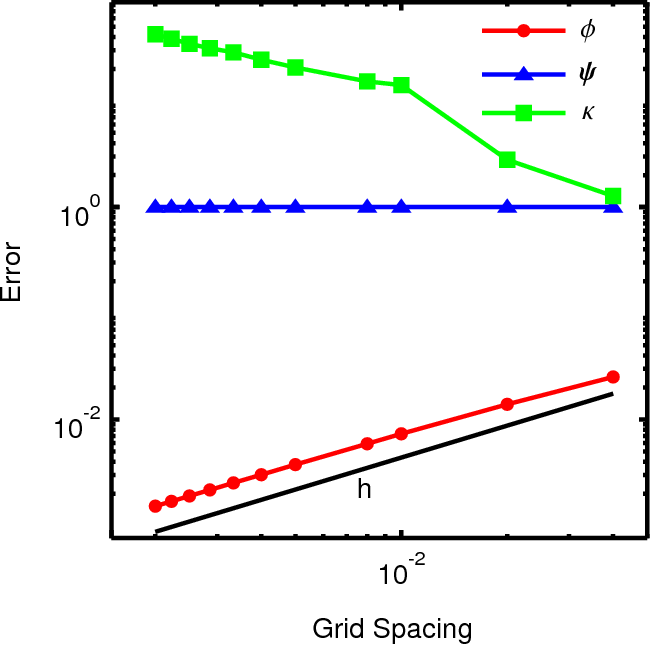}
			}
		\end{center}
		\caption{The error in the entire domain for the ellipse. The $L_\infty$ error does not converge for the gradient and curvature
		fields, as expected by the results of Sec. \ref{sec:expected_errors}.}
 		\label{fig:ellipse_entire_convergence}
	\end{figure}	
			
	\begin{figure}[!ht]
		\begin{center}
			\subfigure[$L_2$ Error]{
				\includegraphics[width=0.475\textwidth]{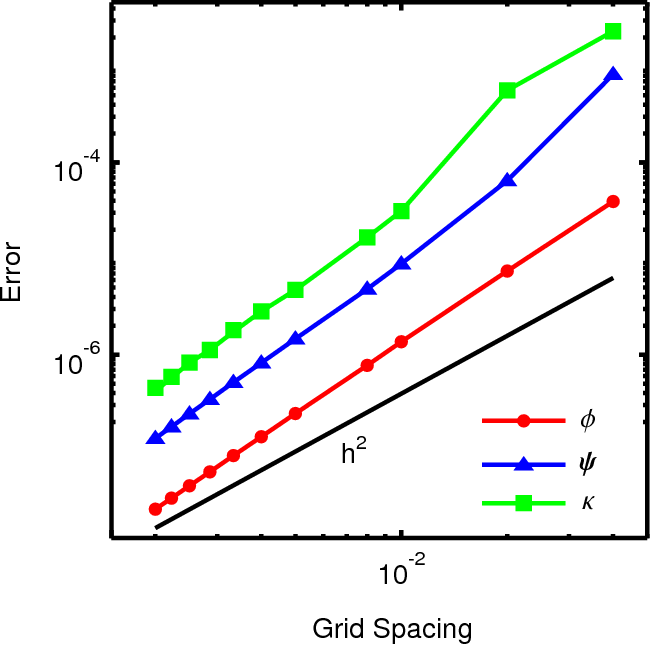}
			}
			\subfigure[$L_\infty$ Error]{
				\includegraphics[width=0.475\textwidth]{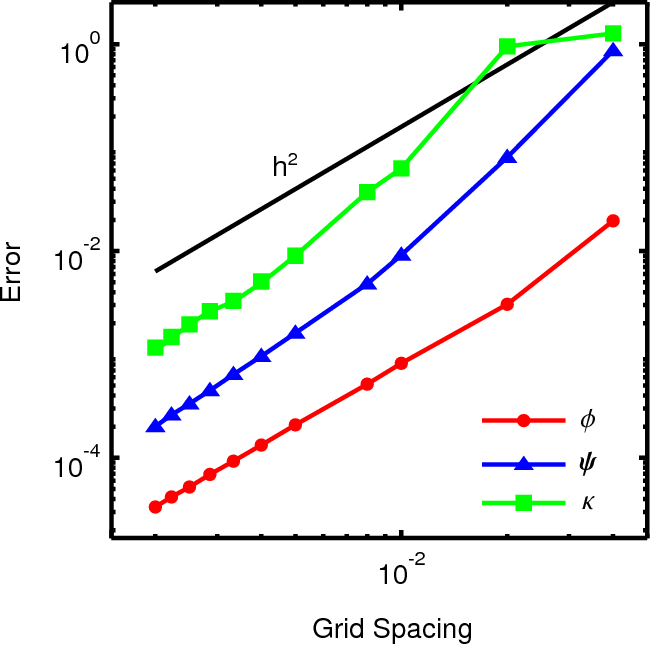}
			}
		\end{center}
		\caption{The error for the ellipse in the region given by $|\phi|\leq 9h$, where $h$ is the grid spacing. In this case
			both the $L_2$ and $L_\infty$ error-norms converge.}		
		\label{fig:ellipse_tube_convergence}
	\end{figure}	
	
\begin{figure}[!ht]
		\begin{center}		
			\subfigure[Level Set: $50^2$]{
				\includegraphics[width=0.475\textwidth]{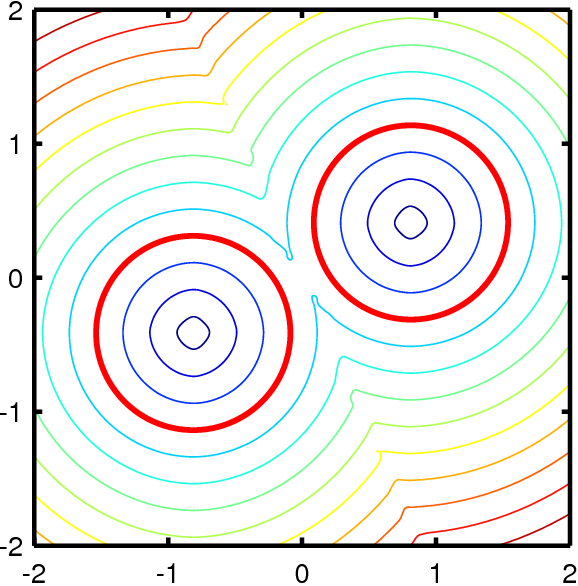}
			}
			\subfigure[Curvature: $50^2$]{
				\includegraphics[width=0.475\textwidth]{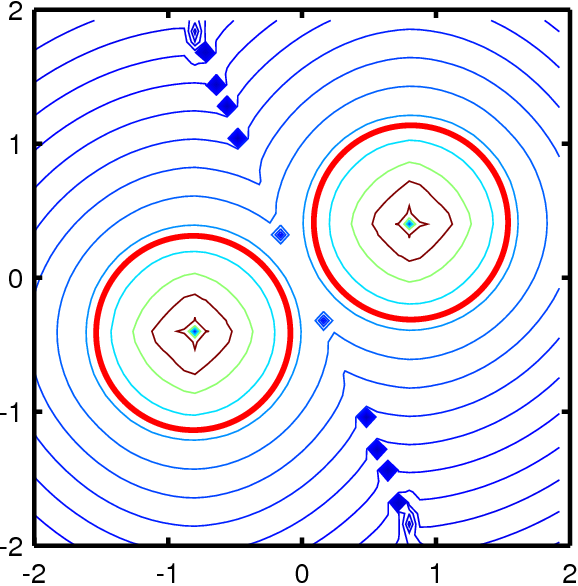}
			}\\
			\subfigure[Level Set: $500^2$]{
				\includegraphics[width=0.475\textwidth]{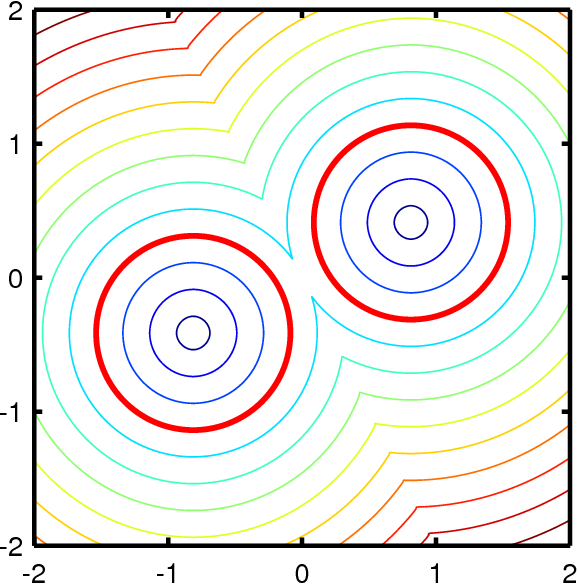}
			}
			\subfigure[Curvature: $500^2$]{
				\includegraphics[width=0.475\textwidth]{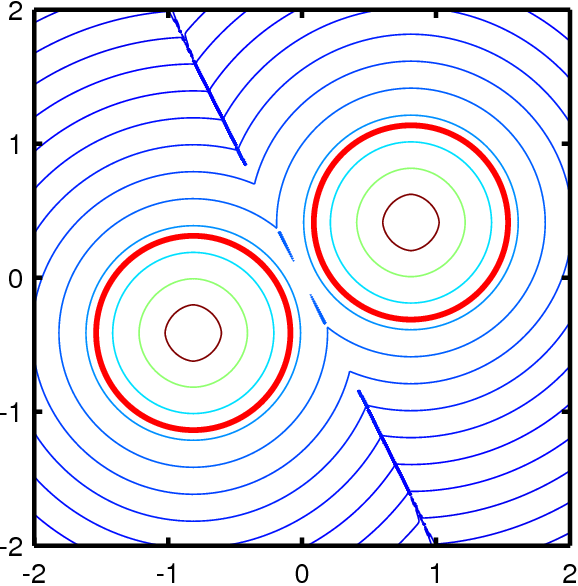}
			}\\
		\end{center}
		\caption{Contours of the level set and curvature for dual circles centered at $(0.8125,0.4125)$ and $(-0.8125,-0.4125)$ and both with radius of 0.75 
			on a coarse and fine grid.}	
		\label{fig:dual_circles}
	\end{figure}

\begin{figure}[!ht]
		\begin{center}
			\subfigure[Level Set: $50^2$]{
				\includegraphics[width=0.475\textwidth]{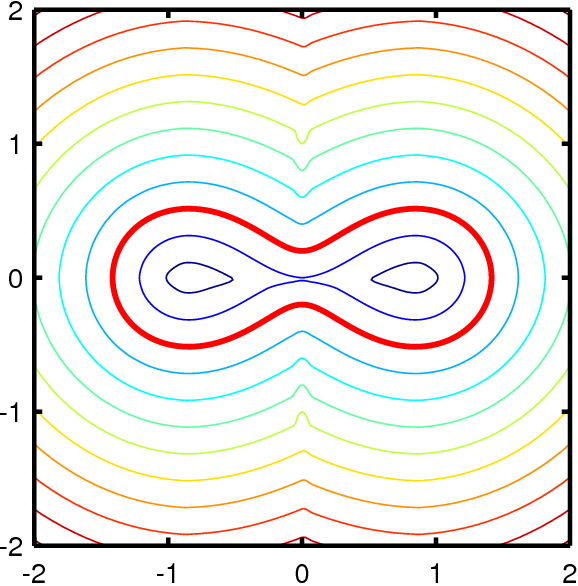}
			}
			\subfigure[Curvature: $50^2$]{
				\includegraphics[width=0.475\textwidth]{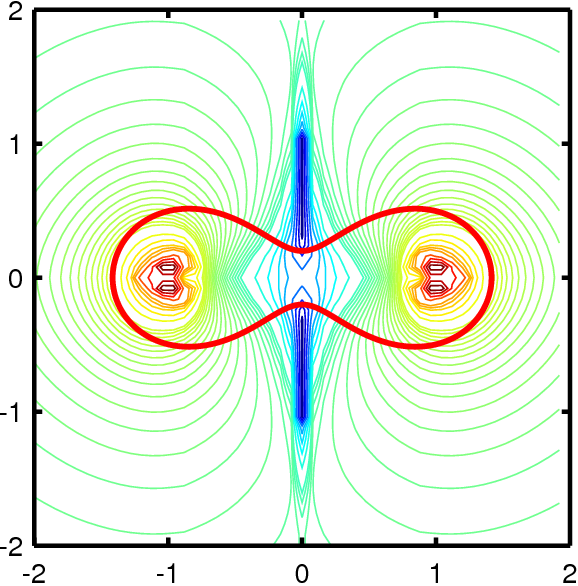}
			} \\			
			\subfigure[Level Set: $500^2$]{
				\includegraphics[width=0.475\textwidth]{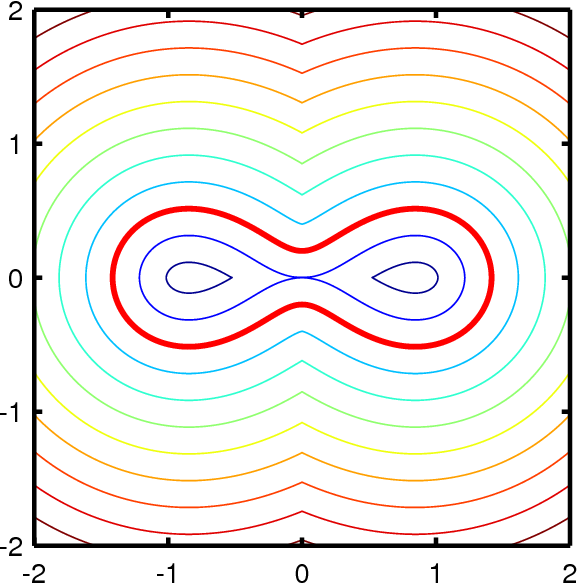}
			}
			\subfigure[Curvature: $500^2$]{
				\includegraphics[width=0.475\textwidth]{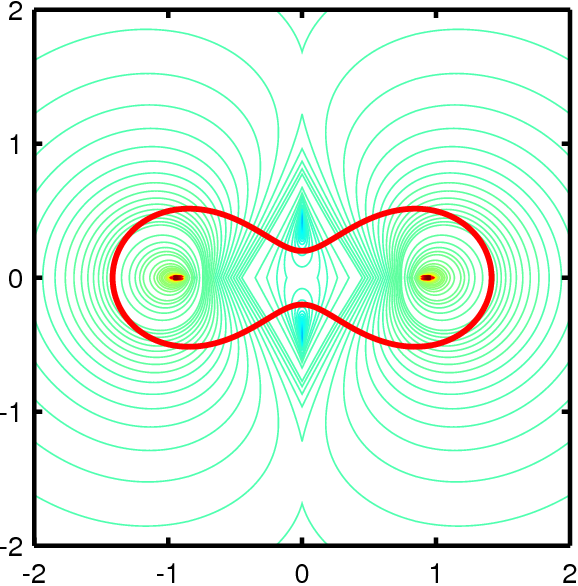}
			} \\
		\end{center}
		\caption{Contours of the level set and curvature for a cassini oval with an initial level set of 
			$\phi_0=\left((x-a)^2+y^2\right)\left((x+a)^2+y^2\right)-b^4$ with $a=0.99$ and $b=1.01$ on a coarse and fine grid.}	
		\label{fig:cassini2d}
	\end{figure}

\begin{figure}[!ht]
		\begin{center}
			\subfigure[Level Set: $100^2$]{
				\includegraphics[width=0.475\textwidth]{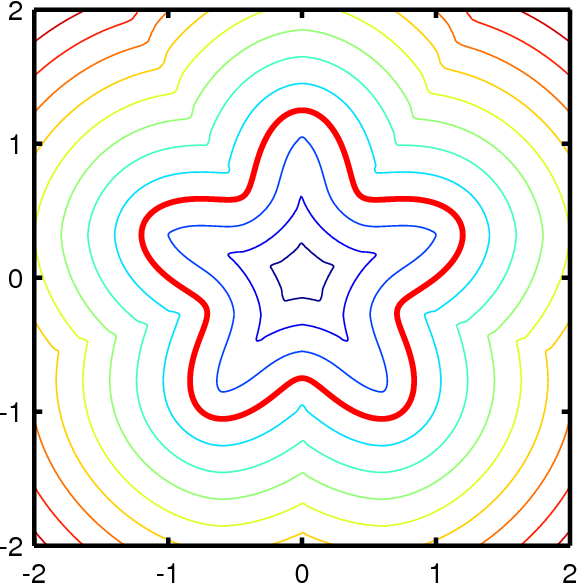}
			}
			\subfigure[Curvature: $100^2$]{
				\includegraphics[width=0.475\textwidth]{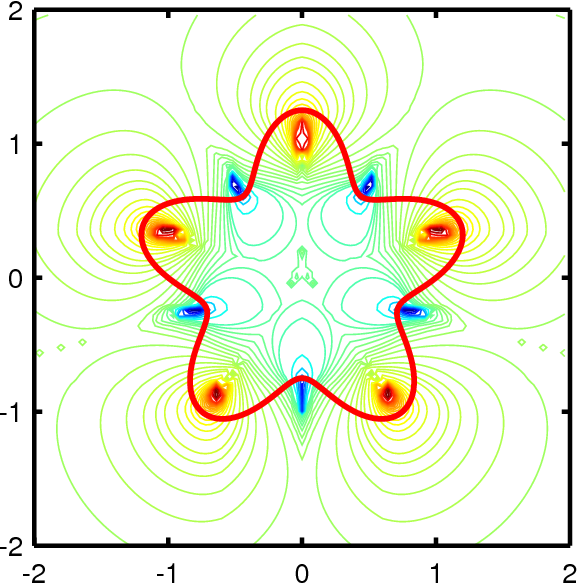}
			} \\
			\subfigure[Level Set: $500^2$]{
				\includegraphics[width=0.475\textwidth]{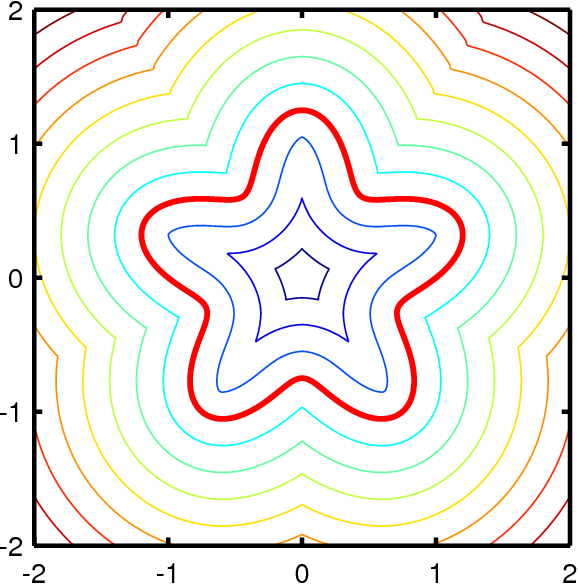}
			}
			\subfigure[Curvature: $500^2$]{
				\includegraphics[width=0.475\textwidth]{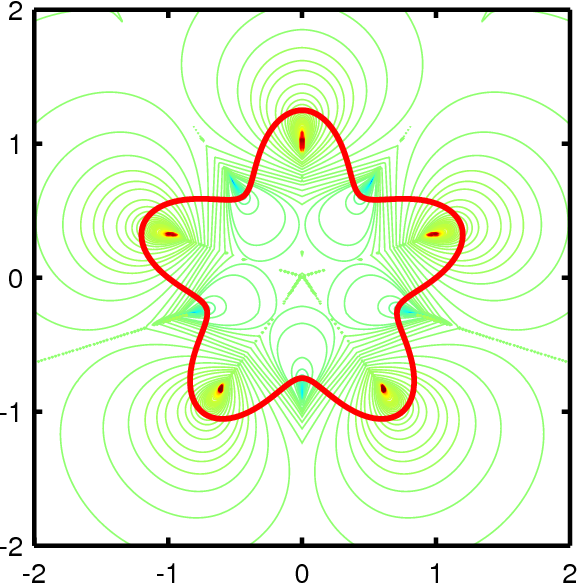}
			} \\
		\end{center}
		\caption{A star interface given by $\phi_0=\sqrt{x^2+y^2}-1+\sin(5\theta)/4$ with $\theta=\rm{ArcTan}(y/x)$  on a coarse and fine grid.}
		\label{fig:star}	
	\end{figure}

	This section concludes by presenting three additional sample interface. In each case a coarse grid of $50^2$ is compared to a fine grid of $500^2$.
	Results are for two circles of radius 0.75 centered at $(0.8125,0.4125)$ and $(-0.8125,-0.4125)$, a Cassini oval with an initial 
	level set of $\phi_0=\left((x-a)^2+y^2\right)+\left((x+a)^2+y^2\right)-b^4$ with $a=0.99$ and $b=1.01$, and a star interface given by
	$\phi_0=\sqrt{x^2+y^2}-1+\sin(5\theta)/4$ with $\theta=\textrm{ArcTan}(y/x)$, see Figs. \ref{fig:dual_circles} to \ref{fig:star}.

	In all three cases level set and curvature fields at the coarse grid qualitatively match those
	at the finer grid. Errors in the curvature field are again demonstrated in regions where level set fronts collide, such as the diagonal in the two circle case of Fig.
	\ref{fig:dual_circles}, the y-axis of the Cassini oval, Fig. \ref{fig:cassini2d}, and along several of the diagonals of the star shape, Fig. \ref{fig:star}. 
	These error diminish as the grid is refined.

\section{Three Dimensional Results}
\label{sec:5.0}
	Sample three-dimensional results are presented here. Due to computational limitations only a limited set of convergence results will be presented
	for three-dimensional surfaces. The domain will be the cube of $[-2,2]^3$ with a uniform grid spacing of $h$ in all directions.
	
\begin{figure}[!ht]
		\begin{center}
			\subfigure[Level Set]{
				\includegraphics[width=0.475\textwidth]{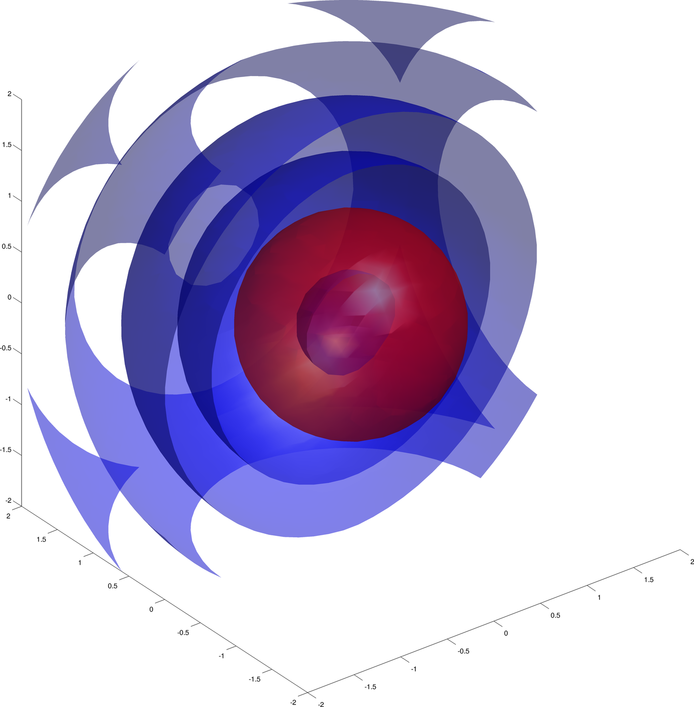}
			}
			\subfigure[Curvature]{
				\includegraphics[width=0.475\textwidth]{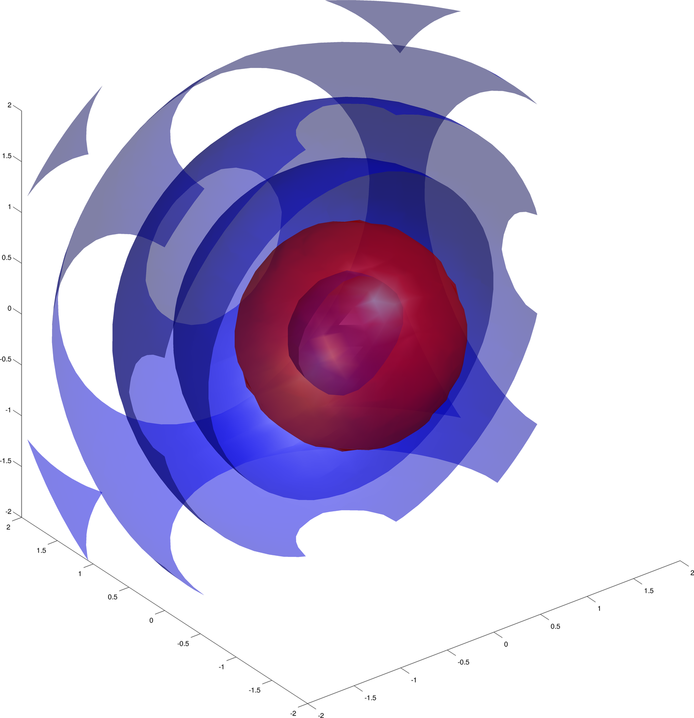}
			}
		\end{center}
		\caption{Sphere of radius one on a $25^3$ grid. The red isosurface represents the value at the interface.}	
		\label{fig:sphere3D_25}
	\end{figure}

\begin{figure}[!ht]
		\begin{center}
			\subfigure[Level Set]{
				\includegraphics[width=0.475\textwidth]{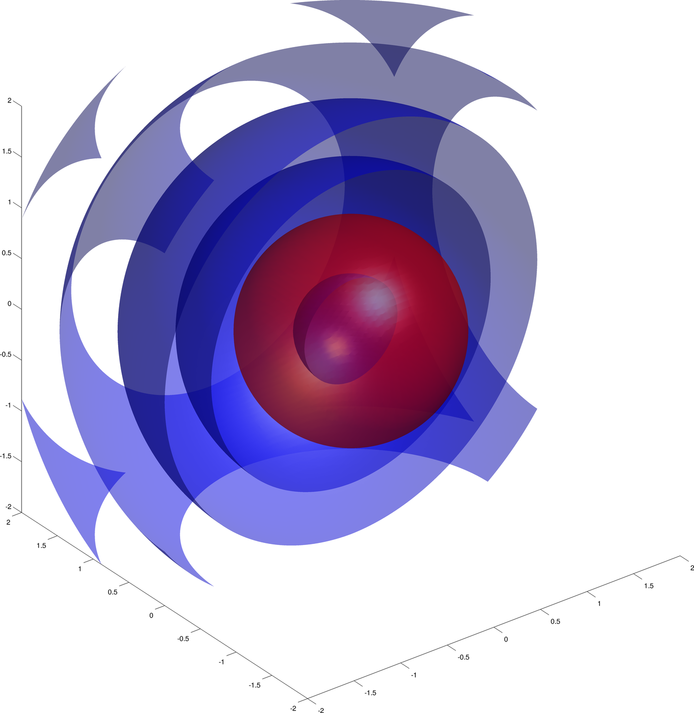}
			}
			\subfigure[Curvature]{
				\includegraphics[width=0.475\textwidth]{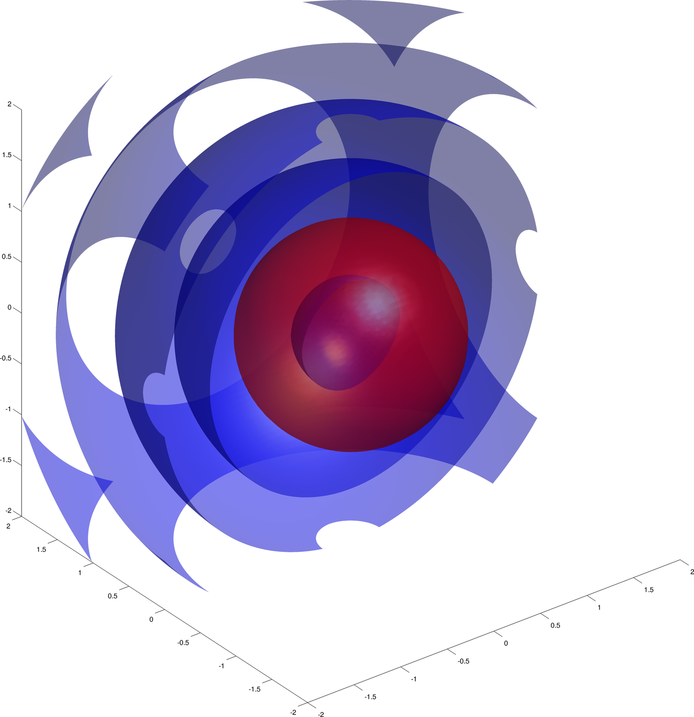}
			}
		\end{center}
		\caption{Sphere of radius one on a $100^3$ grid. The red isosurface represents the value at the interface.}	
		\label{fig:sphere3D_100}		
	\end{figure}

	First consider a spherical surface with a radius of one. The resulting isosurfaces for grids of $25^3$ and $100^3$ are presented in
	Figs. \ref{fig:sphere3D_25} and \ref{fig:sphere3D_100}, respectively. The AFMM on both grids results in a smooth level set and mean curvature field.

\begin{figure}[!ht]
		\begin{center}
			\subfigure[$L_2$ Error]{
				\includegraphics[width=0.475\textwidth]{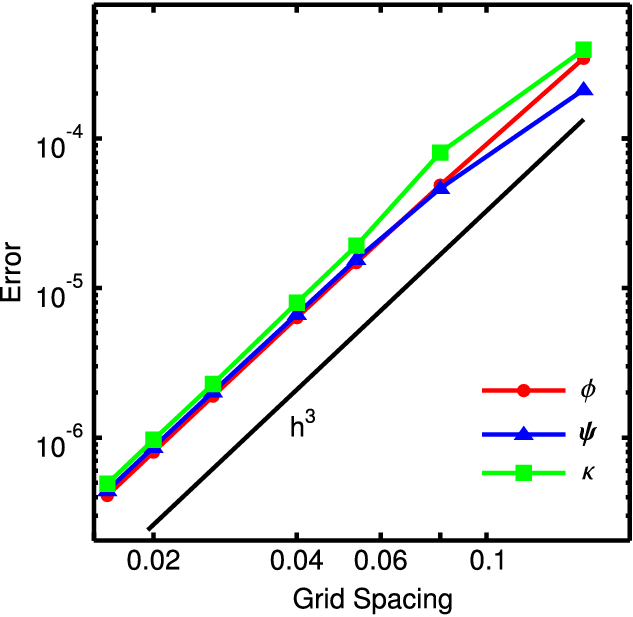}
			}
			\subfigure[$L_\infty$ Error]{
				\includegraphics[width=0.475\textwidth]{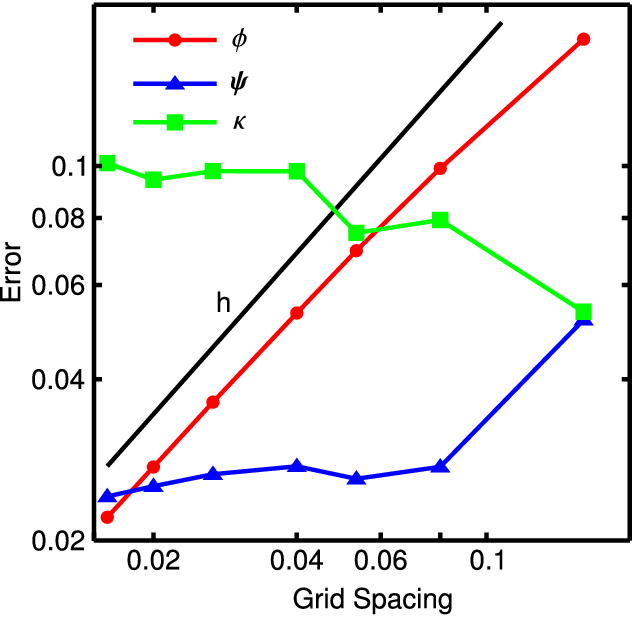}
			}
		\end{center}
		\caption{The error in the entire domain for the spherical surface. The $L_\infty$ error does not converge for the gradient and curvature
		fields, as expected by the results of Sec. \ref{sec:expected_errors}.}
 		\label{fig:sphere3D_entire_convergence}
	\end{figure}

\begin{figure}[!ht]
		\begin{center}
			\subfigure[$L_2$ Error]{
				\includegraphics[width=0.475\textwidth]{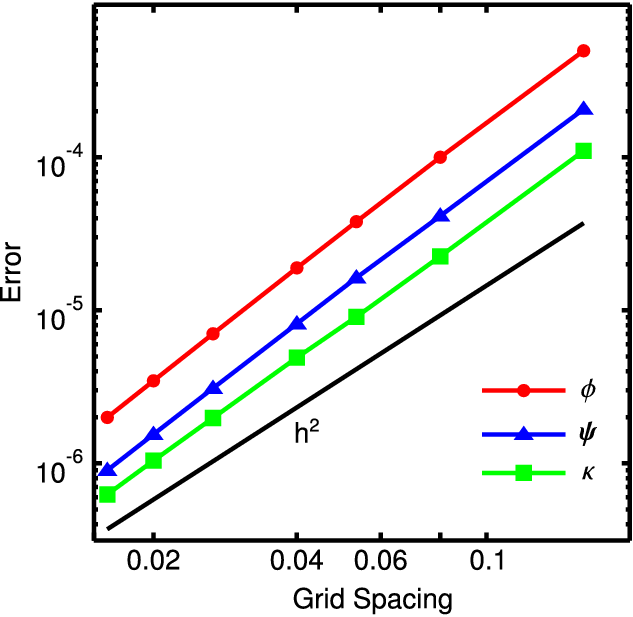}
			}
			\subfigure[$L_\infty$ Error]{
				\includegraphics[width=0.475\textwidth]{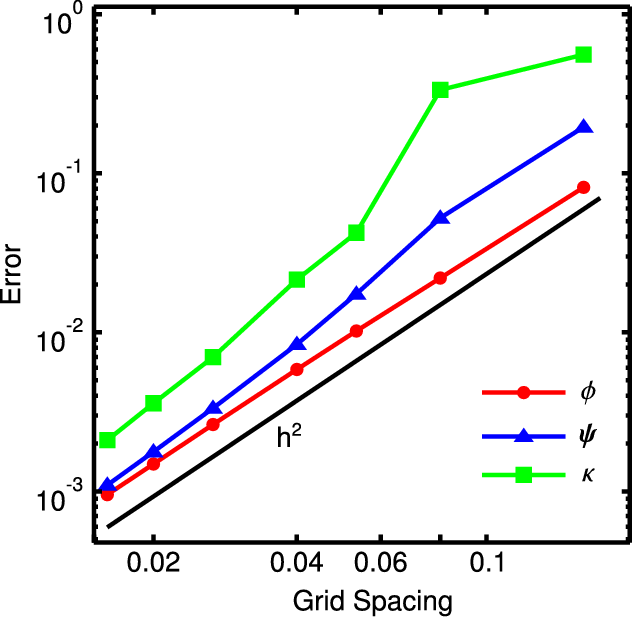}
			}
		\end{center}
		\caption{The error for the spherical surface in the region given by $|\phi|\leq 9h$, where $h$ is the grid spacing. In this case
			both the $L_2$ and $L_\infty$ error-norms converge.}		
		\label{fig:sphere3D_tube_convergence}
	\end{figure}

	Limited convergence results for the spherical interface are presented in Figs. \ref{fig:sphere3D_entire_convergence} and \ref{fig:sphere3D_tube_convergence}. As
	seen in the two-dimensional results the $L_\infty$-norm errors do not converge for the gradient and curvature fields, as expected from Sec. \ref{sec:expected_errors}.
	It appears that in the $L_2$ error for the entire domain converges at third-order. This is most likely due to the additional directions
	that information may travel. The most error will be introduced when only a single neighbor node is available during the calculation of updated values. This is 
	much less likely to occur in three dimensions than in two, resulting in an increase in the accuracy of the scheme.

\begin{figure}[!ht]
		\begin{center}
			\subfigure[Level Set]{
				\includegraphics[width=0.475\textwidth]{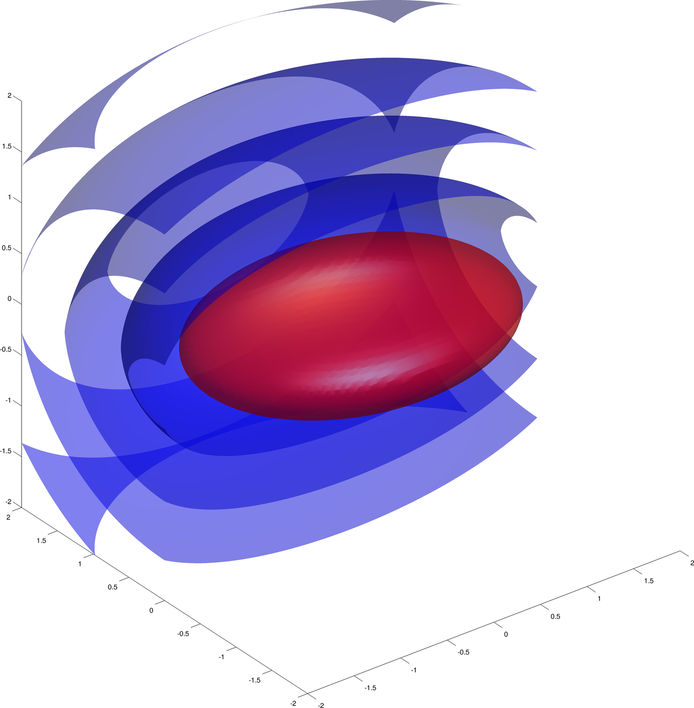}
			}
			\subfigure[Curvature]{
				\includegraphics[width=0.475\textwidth]{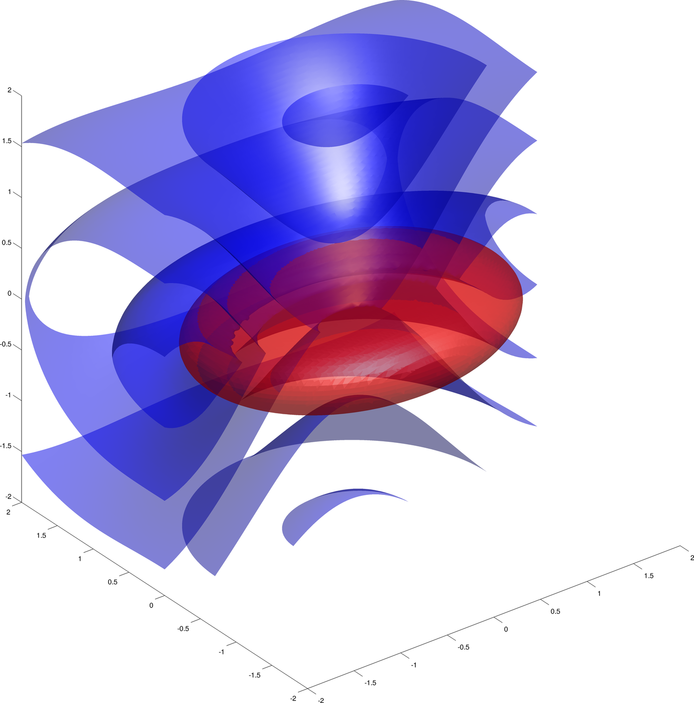}
			}
		\end{center}
		\caption{Ellipsoid given by $\phi_0=\sqrt{(x/1.6)^2+(y/1.2)^2+(z/0.5)^2}-1.0$ on a $100^3$ grid. The red isosurface represents the interface.}
		\label{fig:ellipsoid3D}
	\end{figure}

\begin{figure}[!ht]
		\begin{center}
			\subfigure[Level Set]{
				\includegraphics[width=0.475\textwidth]{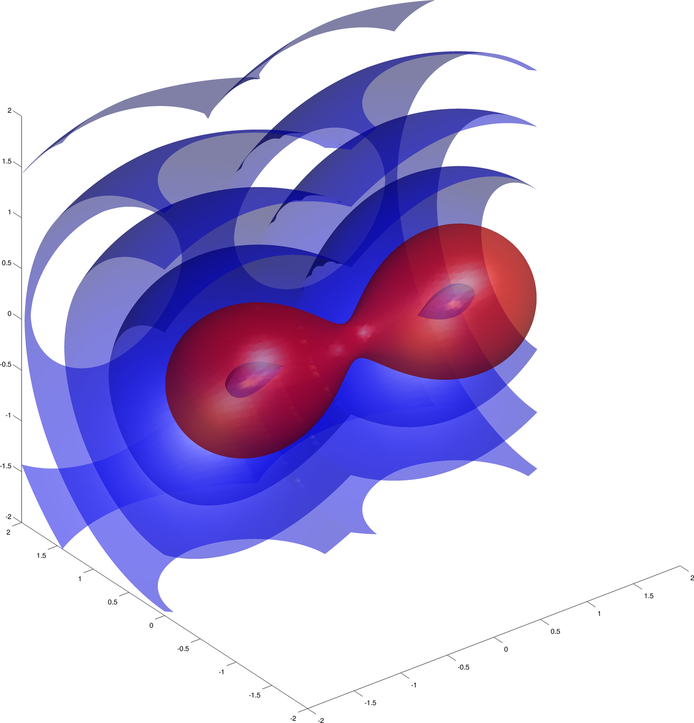}
			}
			\subfigure[Curvature]{
				\includegraphics[width=0.475\textwidth]{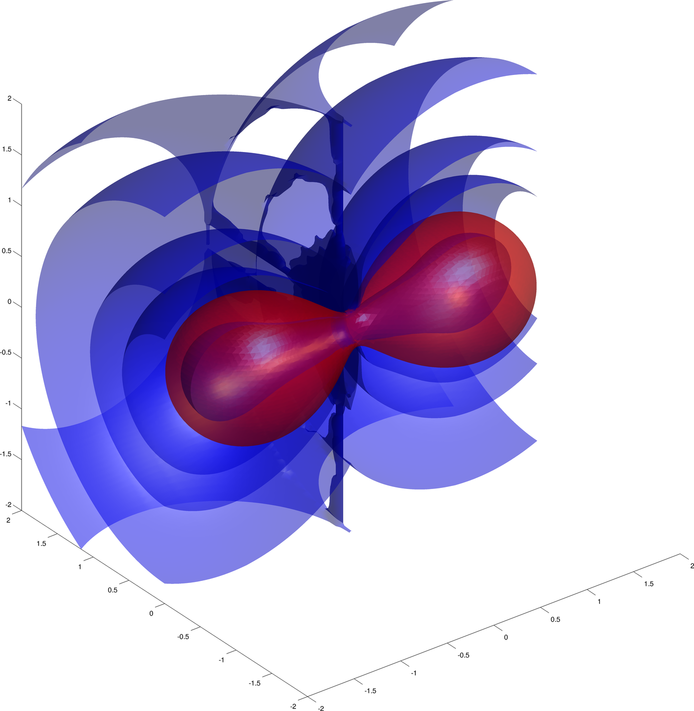}
			}
		\end{center}	
		\caption{Cassini oval given by $\phi_0=\left((x-a)^2+y^2+z^2\right)\left((x+a)^2+y^2+z^2\right)-b^4$ with $a=1.29$ and $b=1.3$ on a $100^3$ grid. The red 
			isosurface represents the interface.}	
		\label{fig:cassini3D}
	\end{figure}

	The result for an ellipsoidal surface given by $\phi_0=\sqrt{(x/1.6)^2+(y/1.2)^2+(z/0.5)^2}-1.0$ using a $100^3$ grid is shown in
	Fig. \ref{fig:ellipsoid3D} while that for a three-dimensional Cassini oval given by 
	$\phi_0=\left((x-a)^2+y^2+z^2\right)\left((x+a)^2+y^2+z^2\right)-b^4$ with $a=1.29$ and $b=1.3$ on a $100^3$ grid is shown in Fig. \ref{fig:cassini3D}.
	In both cases the level set and overall curvature field are smooth. As in the two-dimensional cases 
	slight issues are observed for the Cassini oval in regions where level set contours collide.
	
\section{Conclusion}
\label{sec:6.0}
	In this article the augmented fast marching method for reinitialization of level sets is presented. 
	This work builds upon the fast marching method work of Chopp \cite{Chopp2001} and the gradient augmented level set work of 
	Nave \textit{et. al} \cite{nave2010gradient}. This method increases the accuracy of standard
	level set schemes by calculating the signed distance function and up to second-order derivatives of a general interface. 
	Results show that both the level set and curvature fields are smooth for a wide variety of interfaces and in both two- and three-dimensions. 
	It has also been demonstrated	that the scheme calculates a higher than second-order accurate level set and gradient vector field while the resulting curvature is 
	slightly higher than first-order accurate. This was accomplished using standard first-order upwind derivatives, as in the original fast marching method.
	Unlike partial differential equation based reinitialization schemes a fast marching based method has the advantage of only requiring a single pass to 
	update a level set far from the signed distance function. 

	The additional accuracy has already proven useful in practice. This technique has been used in the modelling of vesicles in general flows 
	\cite{salac2011level, SalacJFM}. In such
	a simulation second-order derivatives of the curvature are required. The use of this reinitialization technique, in conjunction with the gradient-augmented
	level set method, allowed for meaningful simulations to be performed using a single workstation.

\clearpage

\appendix

\section{Two dimensional discretization using only y-direction information}
\label{sec:A}
	The two-dimensional discretization of Eqs. (\ref{eq:2D_pde_phi})-(\ref{eq:2D_pde_phixy}) using only information from the y-direction is 
	\begin{align*}
		\psi^x \psi^x + \psi^y (D_y^{\pm} \phi) &= 1 \\
		\psi^y (D_y^{\pm} \psi^x) & =0 \\
		\psi^y (D_y^{\pm} \psi^x) + \psi^y (D_y^{\pm} \psi^y) & =0 \\
		H^{xx} H^{xx} + H^{xy} H^{xy} + \psi^y (D_y^{\pm} H^{xx}) &=0 \\
		H^{yy} H^{yy} + H^{xy} H^{xy} + \psi^x (D_y^{\pm} H^{xy}) + \psi^y (D_y^{\pm} H^{yy}) &=0 \\
		H^{xx} H^{xy} + H^{yy} H^{xy} + \psi^x (D_y^{\pm} H^{xx}) + \psi^y (D_y^{\pm} H^{xy}) &=0
	\end{align*}		

\section{Three dimensional discretizations}
\label{sec:B}
	The full discretizations for all cases in three dimensions are included here. There are a total of seven cases to consider.
	
\subsection{Information from the x-, y-, and z-directions}
\label{sec:Bxyz}
	\begin{align*}
		\vec{\psi}^x (D_x^\pm \phi)  + \vec{\psi}^y (D_y^\pm \phi) + \vec{\psi}^z (D_z^\pm \phi) &= 1 \\
		\vec{\psi}^x (D_x^\pm \vec{\psi}^x) + \vec{\psi}^y (D_y^\pm \vec{\psi}^x) + \vec{\psi}^z (D_z^\pm \vec{\psi}^x) &= 0 \\
		\vec{\psi}^x (D_x^\pm \vec{\psi}^y) + \vec{\psi}^y (D_y^\pm \vec{\psi}^y) + \vec{\psi}^z (D_z^\pm \vec{\psi}^y) &= 0 \\
		\vec{\psi}^x (D_x^\pm \vec{\psi}^z) + \vec{\psi}^y (D_y^\pm \vec{\psi}^z) + \vec{\psi}^z (D_z^\pm \vec{\psi}^z) &= 0 \\	
		H^{xx} H^{xx} + H^{xy} H^{xy} + H^{xz} H^{xz} + \vec{\psi}^x (D_x^\pm H^{xx}) + \vec{\psi}^y (D_y^\pm H^{xx}) + \vec{\psi}^z (D_z^\pm H^{xx}) &= 0 \\
		H^{xy} H^{xy} + H^{yy} H^{yy} + H^{yz} H^{yz} + \vec{\psi}^x (D_x^\pm H^{yy}) + \vec{\psi}^y (D_y^\pm H^{yy}) + \vec{\psi}^z (D_z^\pm H^{yy}) &= 0 \\
		H^{xz} H^{xz} + H^{yz} H^{yz} + H^{zz} H^{zz} + \vec{\psi}^x (D_x^\pm H^{zz}) + \vec{\psi}^y (D_y^\pm H^{zz}) + \vec{\psi}^z (D_z^\pm H^{zz}) &= 0 \\
		H^{xx} H^{xy} + H^{yy} H^{xy} + H^{xz} H^{yz} + \vec{\psi}^x (D_x^\pm H^{xy}) + \vec{\psi}^y (D_y^\pm H^{xy}) + \vec{\psi}^z (D_z^\pm H^{xy}) &= 0 \\
		H^{xx} H^{xz} + H^{zz} H^{xz} + H^{xy} H^{yz} + \vec{\psi}^x (D_x^\pm H^{xz}) + \vec{\psi}^y (D_y^\pm H^{xz}) +	\vec{\psi}^z (D_z^\pm H^{xz}) &= 0 \\
		H^{yy} H^{yz} + H^{zz} H^{yz} + H^{xz} H^{xy} + \vec{\psi}^x (D_x^\pm H^{yz}) + \vec{\psi}^y (D_y^\pm H^{yz}) + \vec{\psi}^z (D_z^\pm H^{yz}) &= 0
	\end{align*}

\subsection{Information from the x- and y-directions}
\label{sec:Bxy}
	\begin{align*}
		\vec{\psi}^x (D_x^\pm \phi)  + \vec{\psi}^y (D_y^\pm \phi) + \vec{\psi}^z \vec{\psi}^z  &= 1 \\
		\vec{\psi}^x (D_x^\pm \vec{\psi}^x) + \vec{\psi}^y (D_y^\pm \vec{\psi}^x) + \vec{\psi}^z (D_x^\pm \vec{\psi}^z) &= 0 \\
		\vec{\psi}^x (D_x^\pm \vec{\psi}^y) + \vec{\psi}^y (D_y^\pm \vec{\psi}^y) + \vec{\psi}^z (D_y^\pm \vec{\psi}^z) &= 0 \\
		\vec{\psi}^x (D_x^\pm \vec{\psi}^z) + 										\vec{\psi}^y (D_y^\pm \vec{\psi}^z) + &= 0 \\		
		H^{xx} H^{xx} + H^{xy} H^{xy} + H^{xz} H^{xz} + \vec{\psi}^x (D_x^\pm H^{xx}) + \vec{\psi}^y (D_y^\pm H^{xx}) + \vec{\psi}^z (D_x^\pm H^{xz}) &= 0 \\
		H^{xy} H^{xy} + H^{yy} H^{yy} + H^{yz} H^{yz} + \vec{\psi}^x (D_x^\pm H^{yy}) + \vec{\psi}^y (D_y^\pm H^{yy}) + \vec{\psi}^z (D_y^\pm H^{yz}) &= 0 \\
		H^{xz} H^{xz}  + H^{yz} H^{yz} + H^{zz} H^{zz} + \vec{\psi}^x (D_x^\pm H^{zz}) + \vec{\psi}^y (D_y^\pm H^{zz}) &= 0 \\
		H^{xx} H^{xy} + H^{yy} H^{xy} + H^{xz} H^{yz} + \vec{\psi}^x (D_x^\pm H^{xy}) +  \\
				\vec{\psi}^y (D_y^\pm H^{xy}) + \vec{\psi}^z ((D_y^\pm H^{xz}) + (D_x^\pm H^{yz}))/2 &= 0 \\
		H^{xx} H^{xz} + H^{zz} H^{xz} + H^{xy} H^{yz} + \vec{\psi}^x (D_x^\pm H^{xz}) + \vec{\psi}^y (D_y^\pm H^{xz}) + \vec{\psi}^z (D_x^\pm H^{zz}) &= 0 \\
		H^{yy} H^{yz} + H^{zz} H^{yz} + H^{xz} H^{xy} + \vec{\psi}^x (D_x^\pm H^{yz}) + \vec{\psi}^y (D_y^\pm H^{yz}) + \vec{\psi}^z (D_y^\pm H^{zz}) &= 0 
	\end{align*}
\subsection{Information from the x- and z-directions}
\label{sec:Bxz}
	\begin{align*}
		\vec{\psi}^x (D_x^\pm \phi)  + \vec{\psi}^y \vec{\psi}^y + \vec{\psi}^z (D_z^\pm \phi) &= 1 \\
		\vec{\psi}^x (D_x^\pm \vec{\psi}^x) + \vec{\psi}^y (D_x^\pm \vec{\psi}^y) + \vec{\psi}^z (D_z^\pm \vec{\psi}^x) &= 0 \\
		\vec{\psi}^x (D_x^\pm \vec{\psi}^y) + 										\vec{\psi}^z (D_z^\pm \vec{\psi}^y) &= 0 \\
		\vec{\psi}^x (D_x^\pm \vec{\psi}^z) + \vec{\psi}^y (D_z^\pm \vec{\psi}^y) + \vec{\psi}^z (D_z^\pm \vec{\psi}^z) &= 0 \\	
		H^{xx} H^{xx} + H^{xy} H^{xy} + H^{xz} H^{xz} + \vec{\psi}^x (D_x^\pm H^{xx}) + \vec{\psi}^y (D_x^\pm H^{xy}) + \vec{\psi}^z (D_z^\pm H^{xx}) &= 0 \\
		H^{xy} H^{xy} + H^{yy} H^{yy} + H^{yz} H^{yz} + \vec{\psi}^x (D_x^\pm H^{yy}) + \vec{\psi}^z (D_z^\pm H^{yy}) &= 0 \\
		H^{xz} H^{xz}  + H^{yz} H^{yz} + H^{zz} H^{zz} + \vec{\psi}^x (D_x^\pm H^{zz}) + \vec{\psi}^y (D_z^\pm H^{yz}) + \vec{\psi}^z (D_z^\pm H^{zz}) &= 0 \\
		H^{xx} H^{xy} + H^{yy} H^{xy} + H^{xz} H^{yz} + \vec{\psi}^x (D_x^\pm H^{xy}) + \vec{\psi}^y (D_x^\pm H^{yy}) + \vec{\psi}^z (D_z^\pm H^{xy}) &= 0 \\
		H^{xx} H^{xz} + H^{zz} H^{xz} + H^{xy} H^{yz} + \\
				\vec{\psi}^x (D_x^\pm H^{xz}) + \vec{\psi}^y ((D_z^\pm H^{xy}) + (D_x^\pm H^{yz}))/2 + \vec{\psi}^z (D_z^\pm H^{xz}) &= 0 \\
		H^{yy} H^{yz} + H^{zz} H^{yz} + H^{xz} H^{xy} + \vec{\psi}^x (D_x^\pm H^{yz}) + \vec{\psi}^y (D_z^\pm H^{yy}) + \vec{\psi}^z (D_z^\pm H^{yz}) &= 0 
	\end{align*}
\subsection{Information from the y- and z-directions}
\label{sec:Byz}
	\begin{align*}
		\vec{\psi}^x \vec{\psi}^x   + \vec{\psi}^y (D_y^\pm \phi) + \vec{\psi}^z (D_z^\pm \phi) &= 1 \\
											  \vec{\psi}^y (D_y^\pm \vec{\psi}^x) + \vec{\psi}^z (D_z^\pm \vec{\psi}^x) &= 0 \\
		\vec{\psi}^x (D_y^\pm \vec{\psi}^x) + \vec{\psi}^y (D_y^\pm \vec{\psi}^y) + \vec{\psi}^z (D_z^\pm \vec{\psi}^y) &= 0 \\
		\vec{\psi}^x (D_z^\pm \vec{\psi}^x) + \vec{\psi}^y (D_y^\pm \vec{\psi}^z) + \vec{\psi}^z (D_z^\pm \vec{\psi}^z) &= 0 \\	
		H^{xx} H^{xx} + H^{xy} H^{xy} + H^{xz} H^{xz} + \vec{\psi}^y (D_y^\pm H^{xx}) + \vec{\psi}^z (D_z^\pm H^{xx}) &= 0 \\
		H^{xy} H^{xy} + H^{yy} H^{yy} + H^{yz} H^{yz} + \vec{\psi}^x (D_y^\pm H^{xy}) + \vec{\psi}^y (D_y^\pm H^{yy}) + \vec{\psi}^z (D_z^\pm H^{yy}) &= 0 \\
		H^{xz} H^{xz} + H^{yz} H^{yz} + H^{zz} H^{zz} + \vec{\psi}^x (D_z^\pm H^{xz}) + \vec{\psi}^y (D_y^\pm H^{zz}) + \vec{\psi}^z (D_z^\pm H^{zz}) &= 0 \\
		H^{xx} H^{xy} + H^{yy} H^{xy} + H^{xz} H^{yz} + \vec{\psi}^x (D_y^\pm H^{xx}) + \vec{\psi}^y (D_y^\pm H^{xy}) + \vec{\psi}^z (D_z^\pm H^{xy}) &= 0 \\
		H^{xx} H^{xz} + H^{zz} H^{xz} + H^{xy} H^{yz} + \vec{\psi}^x (D_z^\pm H^{xx}) + \vec{\psi}^y (D_y^\pm H^{xz}) + \vec{\psi}^z (D_z^\pm H^{xz}) &= 0 \\
		H^{yy} H^{yz} + H^{zz} H^{yz} + H^{xz} H^{xy} + \\
				\vec{\psi}^x ((D_y^\pm H^{xz}) + (D_z^\pm H^{xy}))/2 + \vec{\psi}^y (D_y^\pm H^{yz}) + \vec{\psi}^z (D_z^\pm H^{yz}) &= 0 
	\end{align*}
\subsection{Information from the x-direction}
\label{sec:Bx}
	\begin{align*}
		(D_x^\pm \phi) \vec{\psi}^x + \vec{\psi}^y \vec{\psi}^y + \vec{\psi}^z \vec{\psi}^z &= 1 \\
		\vec{\psi}^x (D_x^\pm \vec{\psi}^x) + \vec{\psi}^y (D_x^\pm \vec{\psi}^y) + \vec{\psi}^z (D_x^\pm \vec{\psi}^z) &= 0 \\
		\vec{\psi}^x (D_x^\pm \vec{\psi}^y) &= 0 \\
		\vec{\psi}^x (D_x^\pm \vec{\psi}^z) &= 0 \\	
		H^{xx} H^{xx} + H^{xy} H^{xy} + H^{xz} H^{xz} + \vec{\psi}^x (D_x^\pm H^{xx}) + \vec{\psi}^y (D_x^\pm H^{xy}) + \vec{\psi}^z (D_x^\pm H^{xz}) &= 0 \\
		H^{xy} H^{xy} + H^{yy} H^{yy} + H^{yz} H^{yz} + \vec{\psi}^x (D_x^\pm H^{yy}) &= 0 \\
		H^{xz} H^{xz}  + H^{yz} H^{yz} + H^{zz} H^{zz} + \vec{\psi}^x (D_x^\pm H^{zz}) &= 0 \\
		H^{xx} H^{xy} + H^{yy} H^{xy} + H^{xz} H^{yz} + \vec{\psi}^x (D_x^\pm H^{xy}) + \vec{\psi}^y (D_x^\pm H^{yy}) + \vec{\psi}^z (D_x^\pm H^{yz}) &= 0 \\
		H^{xx} H^{xz} + H^{zz} H^{xz} + H^{xy} H^{yz} + \vec{\psi}^x (D_x^\pm H^{xz}) + \vec{\psi}^y (D_x^\pm H^{yz}) + \vec{\psi}^z (D_x^\pm H^{zz}) &= 0 \\
		H^{yy} H^{yz} + H^{zz} H^{yz} + H^{xz} H^{xy} + \vec{\psi}^x (D_x^\pm H^{yz}) &= 0 
	\end{align*}
\subsection{Information from the y-direction}
\label{sec:By}
	\begin{align*}
		\vec{\psi}^x \vec{\psi}^x + (D_y^\pm \phi) \vec{\psi}^y + \vec{\psi}^z \vec{\psi}^z &= 1 \\
		\vec{\psi}^y (D_y^\pm \vec{\psi}^x) &= 0 \\
		\vec{\psi}^x (D_y^\pm \vec{\psi}^x) + \vec{\psi}^y (D_y^\pm \vec{\psi}^y) + \vec{\psi}^z (D_y^\pm \vec{\psi}^z) &= 0 \\
		\vec{\psi}^y (D_y^\pm \vec{\psi}^z) &= 0 \\	
		H^{xx} H^{xx} + H^{xy} H^{xy} + H^{xz} H^{xz} + \vec{\psi}^y (D_y^\pm H^{xx}) &= 0 \\
		H^{xy} H^{xy} + H^{yy} H^{yy} + H^{yz} H^{yz} + \vec{\psi}^x (D_y^\pm H^{xy}) + \vec{\psi}^y (D_y^\pm H^{yy}) + \vec{\psi}^z (D_y^\pm H^{yz}) &= 0 \\
		H^{xz} H^{xz}  + H^{yz} H^{yz} + H^{zz} H^{zz} + \vec{\psi}^y (D_y^\pm H^{zz}) &= 0 \\
		H^{xx} H^{xy} + H^{yy} H^{xy} + H^{xz} H^{yz} + \vec{\psi}^x (D_y^\pm H^{xx}) + \vec{\psi}^y (D_y^\pm H^{xy}) + \vec{\psi}^z (D_y^\pm H^{xz}) &= 0 \\
		H^{xx} H^{xz} + H^{zz} H^{xz} + H^{xy} H^{yz} + \vec{\psi}^y (D_y^\pm H^{xz}) &= 0 \\
		H^{yy} H^{yz} + H^{zz} H^{yz} + H^{xz} H^{xy} + \vec{\psi}^x (D_y^\pm H^{xz}) + \vec{\psi}^y (D_y^\pm H^{yz}) + \vec{\psi}^z (D_y^\pm H^{zz}) &= 0 
	\end{align*}
\subsection{Information from the z-direction}
\label{sec:Bz}
	\begin{align*}
		\vec{\psi}^x \vec{\psi}^x + \vec{\psi}^y \vec{\psi}^y + (D_z^\pm \phi) \vec{\psi}^z &= 1 \\
		\vec{\psi}^z (D_z^\pm \vec{\psi}^x) &= 0 \\
		\vec{\psi}^z (D_z^\pm \vec{\psi}^y) &= 0 \\
		\vec{\psi}^x (D_z^\pm \vec{\psi}^x) + \vec{\psi}^y (D_z^\pm \vec{\psi}^y) + \vec{\psi}^z (D_z^\pm \vec{\psi}^z) &= 0 \\	
		H^{xx} H^{xx} + H^{xy} H^{xy} + H^{xz} H^{xz} + \vec{\psi}^z (D_z^\pm H^{xx}) &= 0 \\
		H^{xy} H^{xy} + H^{yy} H^{yy} + H^{yz} H^{yz} + \vec{\psi}^z (D_z^\pm H^{yy}) &= 0 \\
		H^{xz} H^{xz}  + H^{yz} H^{yz} + H^{zz} H^{zz} + \vec{\psi}^x (D_z^\pm H^{xz}) + \vec{\psi}^y (D_z^\pm H^{yz}) + \vec{\psi}^z (D_z^\pm H^{zz}) &= 0 \\
		H^{xx} H^{xy} + H^{yy} H^{xy} + H^{xz} H^{yz} + \vec{\psi}^z (D_z^\pm H^{xy}) &= 0 \\
		H^{xx} H^{xz} + H^{zz} H^{xz} + H^{xy} H^{yz} + \vec{\psi}^x (D_z^\pm H^{xx}) + \vec{\psi}^y (D_z^\pm H^{xy}) + \vec{\psi}^z (D_z^\pm H^{xz}) &= 0 \\
		H^{yy} H^{yz} + H^{zz} H^{yz} + H^{xz} H^{xy} + \vec{\psi}^x (D_z^\pm H^{xy}) + \vec{\psi}^y (D_z^\pm H^{yy}) + \vec{\psi}^z (D_z^\pm H^{yz}) &= 0 
	\end{align*}


\end{document}